\documentclass[11pt,reqno]{amsart}

\usepackage{graphicx}
\usepackage{url}
\usepackage{amsmath,amssymb,amsthm}

\newcommand{\be}{\begin{equation}}
\newcommand{\ee}{\end{equation}}
\newcommand{\csch}{\mathrm{csch}}

\newcommand{\J}{\mathbf{J}}
\newcommand{\I}{\mathbf{I}}

\newtheorem{prop}{Proposition}

\begin{document}

\title[Theta neuron subject to delayed 
self-feedback]{Theta neuron subject to delayed 
feedback: \\ a prototypical model for self-sustained pulsing}

\author{Carlo R. Laing}
 \email{c.r.laing@massey.ac.nz}
\address{School of Natural and Computational Sciences, 
Massey University,
Private Bag 102-904 
North Shore Mail Centre, 
Auckland 0745,
New Zealand
}%

\author{Bernd Krauskopf}
 \email{B.Krauskopf@auckland.ac.nz}
\address{Department of Mathematics,
University of Auckland,
Private Bag 92019,
Auckland 1142,
New Zealand
}%




\begin{abstract}
We consider a single theta neuron with delayed self-feedback in the form of a Dirac delta function in time. Because the dynamics of a theta neuron on its own can be solved explicitly --- it is either excitable or shows self-pulsations --- we are able to derive algebraic expressions for existence and stability of the periodic solutions that arise in the presence of feedback. These periodic solutions are characterized by one or more equally spaced pulses per delay interval, and there is an increasing amount of multistability with increasing delay time. We present a complete description of where these self-sustained oscillations can be found in parameter space; in particular, we derive explicit expressions for the loci of their saddle-node bifurcations. We conclude that the theta neuron with delayed self-feedback emerges as a prototypical model: it provides an analytical basis for understanding pulsating dynamics observed in other excitable systems subject to delayed self-coupling.  
\end{abstract}








\maketitle

\vspace*{-3mm}
\begin{center}
{\large May 2022}
\end{center}

\section{Introduction}

The theta neuron model, or theta neuron for short, is a mathematical model designed to capture essential features of spiking or bursting neurons \cite{ermkop86}. It takes the form of a differential equation for a single angular variable $\theta(t) \in (-\pi,\pi]$, representing a phase point moving over an attracting periodic orbit; by convention, the theta neuron produces an output spike or pulse when the angle $\theta(t)$ moves through $\pi$. 

We study here a single theta neuron that receives self-feedback of strength $\kappa$ after a time delay $\tau$, where we take the feedback to act instantaneously as modeled by a Dirac delta function in time. The overall model takes the form
\be \frac{d\theta}{dt}=1-\cos{\theta}+(1+\cos{\theta})\left(I+\kappa\sum_{i:-\tau<t_i<0} \delta(t-t_i-\tau)\right),
\label{eq:dth}
\ee
where $\cdots t_{-3}<t_{-2}<t_{-1}<t_0<0$ are the firing times in the past that enter with strength $\kappa$ and after the delay time $\tau$. 

System~\eqref{eq:dth} for $\kappa = 0$ (without self-feedback) is simply 
the equation of a single theta neuron with constant input current $I$. The main advantage of the theta neuron is that its dynamics can be solved explicitly, because its right-hand side is of a particularly simple form. More specifically, it is the (angular) normal form of a saddle-node bifurcation on an invariant-circle (SNIC) bifurcation~\cite{erm96,ermkop86}, which occurs at $I = 0$. One often visualizes the theta neuron as a point moving over the unit circle in the plane, as given by the periodic angle $\theta(t)$. For $I<0$ the system has two equilibria $\theta_\pm=\pm 2\tan^{-1}{(\sqrt{-I})}$, and it is excitable close to the SNIC bifurcation. The point $\theta_-$ is an attractor, while $\theta_+$ is unstable (a saddle on the attracting circle) and acts like a threshold. If an initial condition $\theta(0)$ (caused by some external perturbation) is smaller than $\theta_+$ then $\theta(t)$ relaxes back to $\theta_-$. However, if $\theta(0)$ is (slightly) larger than $\theta_+$ then $\theta(t)$ increases through $\pi$ and, hence, the theta neuron fires before approaching $\theta_-$ (in the absence of any further perturbations). For $I>0$, on the other hand, there are no equilibria and $\theta(t)$ increases monotonically, and the theta neuron fires periodically with period $\pi/\sqrt{I}$. 

This discussion shows that the theta neuron (without feedback) is an excitable system in the parameter range before it bifurcates to producing sustained periodic oscillations. Excitability is a phenomenon that is observed in numerous natural and man-made systems. It is characterized by an all-or-none response (a pulse or spike) to an external input perturbation, depending on whether or not the perturbation exceeds the so-called excitability threshold. An output pulse corresponds to the sudden release of stored energy and it is followed by a refractory period during which energy is replenished and the excitable system is not able to react to a further purturbation~\cite{murray}. Excitability has been found experimentally and in associated mathematical models of neurons and other cells, as well as different types of laser systems; see, for example, \cite{IzhikevichBook,exciteoverview} as entry points to the literature. 

An excitable unit or cell can react to several input perturbations, where the overall strength of all inputs and their timing determine whether an output pulse ensues. In a (neural) network there are necessarily communication delays between cells that need to be taken into account~\cite{pophau06} in order to understand the ensuing dynamics. In the simplest case, a single excitable system receives its own feedback by its output being re-injected after a delay $\tau$. This overall system is able to regenerate its own response after an external input, resulting in feedback-induced self-pulsations whose timing is controlled by the delay time $\tau$. This general mechanism for sustaining pulses in the delay line or delay interval has been demonstrated in a number of laser systems \cite{GarbinNC15, KelleherPRE10, KrauskopfWalker, RomeiraNSR16, TerrienPRA17, Terrien2018OL} and, recently also in an experiment with an actual cell~\cite{wedslo21}. The connection of a neuron to itself is known as an {\em autapse}~\cite{seulee00} and these are known to occur naturally~\cite{tambuh97}. Stable self-generation of spikes or pulses are important, for example, as a means of memory storage \cite{ConnorsN02, ChaudhuriNN16}.

\begin{figure}[t!]
\begin{center}
\includegraphics[width=12cm]{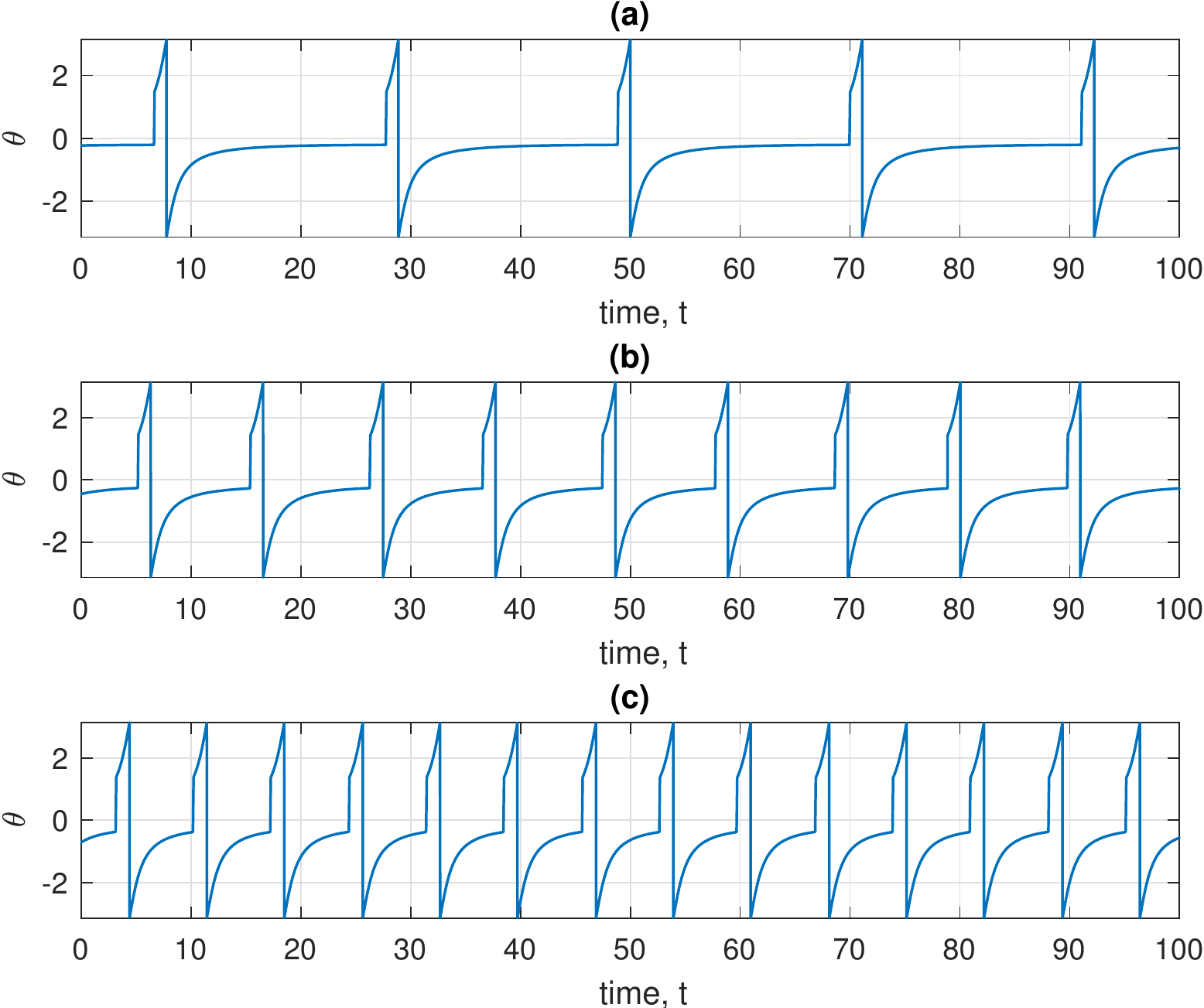}
\caption{Coexisting stable periodic solutions of Eq.~\eqref{eq:dth} for 
$\kappa=1,\tau=20$ and $I=-0.01$. Note that when $\theta$ reaches $\pi$ from below it is
reset to $-\pi$.}
\label{fig:exam}
\end{center}
\end{figure}

We study the theta neuron with self-feedback~\eqref{eq:dth} as a prototypical example of such a system, which allows us to determine the observed dynamics explicitly. The underlying idea is that exact details of the excitability are not important when it comes to identifying underlying principles of self-sustained oscillations in the presence of delayed feedback. Figure~\ref{fig:exam} shows three periodic solutions of Eq.~\eqref{eq:dth} that coexist and are characterized by one, two and three  equidistant pulses being regenerated in the feedback loop over one delay interval. It is the existence and stability of these types of solutions that are the focus of this paper.

A number of other authors have studied oscillators with pulsatile 
interactions~\cite{kliluc18,mirstr90,goeerm02}, exploiting the fact that it is often
possible to explicitly calculate the effects of these interactions.
One example is the leaky integrate-and-fire neuron~\cite{mirstr90}, featuring a reset when the firing threshold is reached, which can be solved explicitly. However, that neuron model has the disadvantage of being nonsmooth and phenomenological, whereas the theta neuron is smooth and has been derived from a more complex Hodgkin-Huxley type model via the technique of phase reduction \cite{erm96,ermkop86}. Note further that~\eqref{eq:dth} is exactly equivalent under the transformation $V=\tan{(\theta/2)}$ to the quadratic integrate-and-fire (QIF) neuron given by the equation
\be
   \frac{dV}{dt}=I+V^2+\kappa\sum_{i:-\tau<t_i<0} \delta(t-t_i-\tau) \label{eq:QIF}
\ee
for the voltage $V$, along with the rule that if $V(t^-)=\infty$ then $V(t^+)=-\infty$~\cite{devmon18,pazmon16}. We will use this equivalence several times below to derive
solutions of the theta neuron model.

We consider here two important subcases of Eq.~\eqref{eq:dth}: the case of $I<0$ when the theta neuron is excitable and has positive (excitatory) self-coupling for $\kappa>0$; and the case of $I>0$ when the theta neuron is intrinsically firing and receives either excitatory (for $\kappa>0$) or inhibitory (for $\kappa<0$) self-coupling. For ease of notation we define the input current $I$ separately as 
\be
   I=\begin{cases} -I_m^2 & I<0 \\ \ \ I_p^2 & I>0 \end{cases} 
\label{eq:I}
\ee
for these two cases. In Sec.~\ref{sec:excit} we consider the case of $-I_m^2 = I < 0$ and $\kappa>0$ of the excitable theta neuron with excitatory self-feedback; we derive the existence of periodic orbits, determine their stability analytically, and give an explicit expression for the saddle-node bifurcations of these orbits. Section~\ref{sec:pos} concerns the case $I_p^2  = I > 0$ of a theta neuron firing periodically even without any feedback. We derive existence and stability of periodic orbits for $\kappa>0$ and then show that these are related to those for $\kappa<0$ via a simple geometric transformation; again, curves of saddle-node bifurcations are given explicitly. In Sec.~\ref{sec:smooth} we consider the theta neuron with a feedback term that is smooth, rather than a delta function. The resulting system needs to be studied numerically, and we find that for both $I<0$ and $I>0$ the dynamics for excitatory coupling with $\kappa>0$ are qualitatively as those of~\eqref{eq:dth}. For $I>0$ and $\kappa<0$, on the other hand, we find additional bifurcations that may lead to chaotic behaviour; however, it can be argued that such a model may not be an appropriate description of a neuron operating in this regime of an oscillating neuron with inhibitory self-feedback. We conclude with a discussion in Sec.~\ref{sec:disc}.

\section{Excitable theta neuron for negative current}
\label{sec:excit}

We first consider the case $-I_m^2 = I < 0$. In order to have nontrivial solutions we consider excitatory coupling only.  Our starting point is the derivation of the solution of~\eqref{eq:dth} for $\kappa=0$, using its equivalence to~\eqref{eq:QIF}. The solution of~\eqref{eq:QIF} for $-I_m<V(0)<I_m$ is
\be
   V(t)=-I_m\tanh{\left(I_mt-\tanh^{-1}{[V(0)/I_m]}\right)}.
\ee
Using the transformation $V=\tan{(\theta/2)}$ it follows that if $-I_m<\tan{(\theta(0)/2)}<I_m$, i.e.~the neuron's state is initially between the
fixed points, then the solution of~\eqref{eq:dth} is
\be 
   \theta(t)=2\tan^{-1}{\left[-I_m\tanh{\left(I_mt-\tanh^{-1}{\left[\frac{\tan{(\frac{\theta(0)}{2})}}{I_m}\right]}\right)}\right]},
\ee
and $\theta$ approaches the lower fixed point from above as $t\to\infty$. Conversely, if $I_m<\tan{(\theta(0)/2)}$, i.e.~the neuron's state is initially above the upper fixed point, then the solution of~\eqref{eq:dth} for $\kappa=0$ is
\be 
    \theta(t)=2\tan^{-1}{\left[-I_m\coth{\left(I_mt-\coth^{-1}{\left[\frac{\tan{(\frac{\theta(0)}{2})}}{I_m}\right]}\right)}\right]}.
\label{eq:sol}
\ee
This means that $\theta$ increases through $\pi$ (the neuron fires)
and $\theta$ then approaches the lower fixed point from below as
$t\to\infty$. While these expressions are general, we will often use them in cases where they simplify. 

Suppose now that the delta function acts at time $t$ and $\kappa\neq 0$. From~\eqref{eq:QIF} and the relationship $V=\tan{(\theta/2)}$ we see that this has the effect of instantaneously changing $V$ by the rule $V(t^+)=V(t^-)+\kappa$, i.e.~$\theta$ is changed by
\[ \tan{\left(\frac{\theta(t^+)}{2}\right)}=\tan{\left(\frac{\theta(t^-)}{2}\right)}+\kappa,
\]
where superscripts indicate the state immediately before ($-$) or after ($+$) the action of the delta function.

\subsection{Existence of periodic solutions for negative current}

We now derive the existence of ocillations such as those in Fig.~\ref{fig:exam}. Suppose the neuron fires at time 0 (i.e.~$\theta(0)=\pi$) and there are an additional $n$ past firing times in the interval $(-\tau,0)$ (here $n$ may be zero); in other words, there are $n+1$ spikes in the delay interval $(-\tau,0]$.
Assume that the past firing times are evenly spaced with time $T$ between them, meaning that $T$ is the period of oscillation. It follows that
\[
   \frac{\tau}{n+1}<T<\frac{\tau}{n}.
\]
The time until the neuron feels the next delta function is $\tau-nT$, which
is less than $T$. Since $\theta(0)=\pi$ we have from~\eqref{eq:sol} that
\[
   \theta(\tau-nT)=2\tan^{-1}{\left[-I_m\coth{\left[I_m(\tau-nT)\right]}\right]}.
\]
The delta function moves $\theta$ to the new value:
\be \theta(\tau-nT^+)=2\tan^{-1}{\left[\tan{\left(\frac{\theta(\tau-nT^-)}{2}\right)}+\kappa\right]}.
\label{eq:thp}
\ee
Assuming that $\tan{[\theta(\tau-nT^+)/2]}>I_m$, i.e.~we are above the upper fixed point, the phase $\theta$ continues to increase and the neuron fires after a further time $\Delta$. 
Again using~\eqref{eq:sol}, this happens when
\[ \pi=2\tan^{-1}{\left[-I_m\coth{\left(I_m\Delta-\coth^{-1}{\left[\frac{\tan{(\frac{\theta(\tau-nT^+)}{2})}}{I_m}\right]}\right)}\right]},
\]
or, equivalently, when
\be
I_m\Delta-\coth^{-1}{\left[\frac{\tan{(\frac{\theta(\tau-nT^+)}{2})}}{I_m}\right]} = 0. 
\label{eq:pi}
\ee
Substituting~\eqref{eq:thp} into~\eqref{eq:pi} we get
\[
   \Delta=\frac{1}{I_m}\coth^{-1} \left[\kappa/I_m-\coth{[I_m(\tau-nT)]}\right].
\]
We know that the amount of time the neuron waits before the delta function acts, ($\tau-nT$) plus $\Delta$, has to equal $T$, i.e.~$\Delta+\tau-nT=T$. Thus we have $\Delta=(n+1)T-\tau$ and finally
\be (n+1)T=\tau+\frac{1}{I_m}\coth^{-1}\left[\frac{\kappa}{I_m}-\coth{[I_m(\tau-nT)]}\right]. 
\label{eq:exist}
\ee
Recall that this expression is valid only for $\frac{\tau}{n+1}<T<\frac{\tau}{n}$. For $n=0$~\eqref{eq:exist} gives the period $T$ explicitly in terms of the other parameters. The periodic solutions shown in Fig.~\ref{fig:exam} with one up to three equidistant spikes per delay interval of length $\tau = 20$ in panels (a), (b) and (c) correspond to $n=0,1$ and 2, respectively.

Note from~\eqref{eq:exist} that we can rescale parameters to scale one of $I_m,\kappa$ and $\tau$ to 1. From now on we set $I_m=1$ and we investigate the effect of varying $\tau$ and $\kappa$, the control parameters of the feedback term. For convenience we rewrite~\eqref{eq:exist} as
\be
   \coth{[(n+1)T-\tau]}=\kappa+\coth{[nT-\tau]}.
\label{eq:existB}
\ee
Since $nT-\tau<0$, we have $\coth{[nT-\tau]}<-1$, and since $(n+1)T-\tau>0$, we have $\coth{[(n+1)T-\tau]}>1$. Thus in order to satisfy~\eqref{eq:existB} we must have $\kappa>2$. Another way to see this is that $V$ has the values $\pm 1$ at the equilibria, and $\kappa$ must be larger than the ``gap'' between them, which is of size~2. 

\subsection{Stability of periodic solutions for negative current}
\label{sec:stab1}

We now consider the stability of solutions satisfying~\eqref{eq:existB}. As above, assume that the neuron has just fired at time $t_0$ and there are $n$ past firing times $t_{-1}$ down to $t_{-n}$ in $(t_0-\tau,t_0)$. 
We wait a time $\tau-(t_0-t_{-n})$ until the delta function acts, which maps $\theta$ from
\[ \theta(\tau-t_0+t_{-n}^-)=-2\tan^{-1}{\left[\coth{\left(\tau-t_0+t_{-n}\right)}\right]}
\]
to
\[ \theta(\tau-t_0+t_{-n}^+)=2\tan^{-1}{\left[\tan{\left(\frac{\theta(\tau-t_0+t_{-n}^-)}{2}\right)}+\kappa\right]}.
\]
We then wait a time $\Delta$ until the neuron fires at time $t_1$, where
\[ \Delta-\coth^{-1}{\left[\tan{\left(\frac{\theta(\tau-t_0+t_{-n}^+)}{2}\right)}\right]}=0.
\]
Thus
\[ t_1=t_0+(\tau-t_0+t_{-n})+\Delta=\tau+t_{-n}+\coth^{-1}\left[\kappa-\coth{(\tau-t_0+t_{-n})}\right].
\]
This gives $t_1$ in terms of previous firing times, but the calculation is general. Hence, we have
\be t_i=\tau+t_{i-n-1}+\coth^{-1}\left[\kappa-\coth{(\tau-t_{i-1}+t_{i-n-1})}\right], 
\label{eq:firing}
\ee
that is, a map giving $t_i$ in terms of previous firing times, which we write as $F(t_i,t_{i-1},\dots, t_{i-n-1})=0$. Note that assuming $t_i=iT$ we recover~\eqref{eq:existB}. 

To calculate stability we perturb $t_i$ to $t_i+\eta_i$ about a periodic solution. 
To first order we have
\[
   \frac{\partial F}{\partial t_i}\eta_i+\frac{\partial F}{\partial t_{i-1}}\eta_{i-1}+\cdots+\frac{\partial F}{\partial t_{i-n-1}}\eta_{i-n-1}=0,
\]
which we write as
\be
   \begin{pmatrix} a_{1,i} \\ a_{2,i} \\ a_{3,i} \\ \vdots \\ a_{n+1,i} \end{pmatrix}
=\begin{pmatrix} 0 & 1 & 0 & 0 & \cdots & 0 \\ 0 & 0 & 1 & 0 & \cdots & 0 \\ \vdots & \vdots & \vdots & \vdots & \vdots & \vdots \\ 0 & 0 & \cdots & \cdots & 0 & 1 \\ F_{n+1} & F_n & \cdots & \cdots & F_2 & F_1\end{pmatrix} \begin{pmatrix} a_{1,i-1} \\ a_{2,i-1} \\ a_{3,i-1} \\ \vdots \\ a_{n+1,i-1} \end{pmatrix} 
\label{eq:a}
\ee
where
\[
   F_k=-\frac{\partial F/\partial t_{i-k}}{\partial F/\partial t_{i}}.
\]
Equation~\eqref{eq:a} is of the form ${\bf a}_i= \J\, {\bf a}_{i-1}$ where ${\bf a}_i \in\mathbb{R}^{n+1}$, so the growth or decay of its solutions (and thus the (in)stability of the periodic solution) is determined by the magnitude of the eigenvalues of the matrix $\J$.

Differentiating~\eqref{eq:firing} we find that $\partial F/\partial t_{i}=-1$, $F_k=0$
for $2\leq k\leq n$, and
$F_{n+1}=1-F_1$. Now
\[
   F_1=\frac{\partial F}{\partial t_{i-1}}=\frac{1-\coth^2{(\tau+t_{i-n-1}-t_{i-1})}}{1-[\kappa-\coth{(\tau+t_{i-n-1}-t_{i-1})}]^2}.
\]
For a periodic solution $t_{i-1}-t_{i-n-1}=nT$, so we can express $F_1$ in term of $T$ and the $n$. Since this quantity determines the stability of periodic solutions, we denote it $\gamma$ from now on for notational convenience, given as
\be
   \gamma = F_1=\frac{\coth^2{(\tau-nT)}-1}{[\kappa-\coth{(\tau-nT)}]^2-1}.
\label{eq:gamma}
\ee
We see that $\gamma$ is always positive: the numerator is clearly positive, and from~\eqref{eq:existB} we have that $\kappa+\coth{(nT-\tau)}=\kappa-\coth{(\tau-nT)}>1$, so the denominator is also positive.
The parameter $\gamma$ is a function of the other parameters in the model: the feedback parameters $\kappa$ and $\tau$, the integer $n$ which specifies the form of the solution we consider (with $n+1$ equidistant spikes), and also the period $T$ which, in turn, is determined by the values of $\kappa$, $\tau$ and $n$.
 
The matrix $\J$ of Eq.~\eqref{eq:a} is thus
\be
   \J=\begin{pmatrix} 0 & 1 & 0 & 0 & \cdots & 0 \\ 0 & 0 & 1 & 0 & \cdots & 0 \\ \vdots & \vdots & \vdots & \vdots & \vdots & \vdots \\ 0 & 0 & \cdots & \cdots & 0 & 1 \\ 1-\gamma & 0 & \cdots & \cdots & 0 & \gamma\end{pmatrix}.
\label{eq:J}
\ee
To find out more about the eigenvalues of $\J$ we consider the determinant of $\J-\lambda \I$, which can be expanded down the first column to obtain 
\[
|\J-\lambda \I|=(-\lambda)^{n+1}+\gamma(-\lambda)^n+(\gamma-1)(-1)^{n+1}.
\]
For even $n$ we have $|\J-\lambda \I| =-\lambda^{n+1}+\gamma\lambda^n+1-\gamma$, and  for odd $n$ we have $|\J-\lambda \I|=\lambda^{n+1}-\gamma\lambda^n-1+\gamma$. In both cases the condition $|\J-\lambda \I|=0$ can be written as $g(\lambda)=0$ where 
\be
g(\lambda)=\lambda^{n+1}-\gamma\lambda^n-1+\gamma.
\label{eq:defg}
\ee
It is the roots of this polynomial that determine the stability of a periodic solution.

We see from \eqref{eq:defg} that $\lambda=1$ is always a root of $g$, reflecting the invariance with respect to uniform time translation of the original system. In particular, $g(\lambda)=(\lambda-1)h(\lambda)$
where
\begin{eqnarray} h(\lambda) & = & \lambda^n+(1-\gamma)\lambda^{n-1}+(1-\gamma)\lambda^{n-2}+\cdots+(1-\gamma)\lambda+1-\gamma \nonumber \\
   & = & \lambda^n+(1-\gamma)\sum_{i=0}^{n-1}\lambda^i \nonumber \\
   & = & \lambda^n+(1-\gamma)\frac{1-\lambda^n}{1-\lambda}
\end{eqnarray}
Note that $h(\lambda) \equiv 1$ for $n=0$. The properties of the roots of $h(\lambda)$ were studied in~\cite{kliluc15}; see also~\cite{kliluc15A,kliluc18}. These authors studied a phase oscillator
with delayed pulsatile self-feedback but in a model for which, in the absence of feedback,
$d\theta/dt=1$. Their results are not equivalent to ours since in the absence of feedback our model
approaches a fixed point. We now collect some useful results.

\begin{prop}[Properties of the eigenvalues at a periodic solution]
\verb+ + \\[-0mm]  
\label{prop:g}

\begin{enumerate}
\item If $\gamma=0$ the roots of $g(\lambda)$ are the $(n+1)$st roots of unity, which all have magnitude 1.
\item When $\gamma=1$ we have $g(\lambda)=\lambda^n(\lambda-1)$, i.e.~the eigenvalues of $\J$ are 1 and 0 (with multiplicity $n$) and the orbit is superstable. 
\item Since $0<\gamma$ the periodic orbits never undergo Neimark-Sacker bifurcations; this is unlike the orbits in~\cite{kliluc15,kliluc15A,kliluc18}. Neither do they undergo period-doubling bifurcations.
\item A root of $g(\lambda)$ leaves the unit circle transversely as $\gamma$ increases through  $(n+1)/n$.
\item For $0<\gamma<(n+1)/n$ all eigenvalues have magnitude less than one and the only instability of a periodic orbit occurs when $h(1)=0$, i.e.~when
$\gamma=(n+1)/n$ (with $n>0$).
\item $dT/d\tau=0$ at the point of superstability (where $\gamma=1$). 
\end{enumerate}

\end{prop}

\noindent
{\bf Proofs of statements:} \\[-0mm]

\begin{enumerate}
\item If $\gamma=0$ then $g(\lambda)=\lambda^{n+1}-1$ and the result follows.
\item Substituting $\gamma=1$ into~\eqref{eq:defg} shows that 0 is an eigenvalue with multiplicity $n$. The orbit is called superstable because --- to linear order --- any perturbations to the firing times immediately decay to zero in one iteration.
\item Substituting $\lambda=e^{i\phi}$ into~\eqref{eq:defg} and setting $g(\lambda)$ equal to zero gives
\be
   e^{i(n+1)\phi}+\gamma-1=\gamma e^{in\phi}. 
\label{eq:proof1}
\ee
Taking absolute values gives $|e^{i(n+1)\phi}+\gamma-1|=|\gamma|$. For $\gamma\neq 1$ this
gives $e^{i(n+1)\phi}=1$, i.e.~$\phi=2\pi k/(n+1)$ for some $k\in\mathbb{Z}$. Substituting this
into~\eqref{eq:proof1} we get $\gamma=\gamma e^{2\pi i kn/(n+1)}$, and if $\gamma\neq 0$ then
we must have $k=(n+1)q$ for some $q\in\mathbb{Z}$, i.e.~$\phi$ is a multiple of $2\pi$ and $\lambda=1$. So if $\gamma\not\in\{0,1\}$ then $\lambda=1$ is the only root of $g$ with magnitude 1. The case $\gamma=0$ is ruled out and the case $\gamma=1$ is covered by (ii).
\item 
Differentiating $g(\lambda)=0$ with respect to $\gamma$ we have
\begin{eqnarray}
 &  (n+1)\lambda^n\frac{d\lambda}{d\gamma}-\left(\lambda^n+\gamma n\lambda^{n-1}\frac{d\lambda}{d\gamma}\right)+1=0 \nonumber \\
& \Rightarrow \quad \frac{d\lambda}{d\gamma}=\frac{\lambda^n-1}{\lambda^{n-1}[\lambda(n+1)-\gamma n]}.
\end{eqnarray}
Now if $\gamma\neq (n+1)/n$ (i.e.~away from this instability), we have
$\frac{d\lambda}{d\gamma}|_{\lambda=1}=0$, reflecting that the root $\lambda=1$ of $g$ is always present.
However, if $\gamma=(n+1)/n$ we have
\[ \frac{d\lambda}{d\gamma} = \frac{\lambda^n-1}{(n+1)(\lambda^n-\lambda^{n-1})},
\]
which is undefined at $\lambda=1$.
With L'Hopital's rule we have
\[
   \lim_{\lambda\to 1}\frac{d\lambda}{d\gamma}=\lim_{\lambda\to 1}\frac{n\lambda^{n-1}}{(n+1)[n\lambda^{n-1}-(n-1)\lambda^{n-2}]}=\frac{n}{n+1}>0.
\]
This is the speed at which the root leaves the unit circle.
\item The only instability occurs when 
$h(1)=0$, i.e.~when $\gamma=(n+1)/n$. When $\gamma=1$ the orbit is stable, so it must be stable for $0<\gamma<(n+1)/n$.
\item 
For superstability we
have $\gamma=1$, which (from the definition of $\gamma$ in~\eqref{eq:gamma}) implies that 
\be
   \kappa=2\coth{(\tau-nT)}.  
\label{eq:kap2}
\ee
The equation for the existence of a periodic orbit is~\eqref{eq:existB}, and 
differentiating it with respect to $\tau$ we get
\be \csch^2{[(n+1)T-\tau]}\left((n+1)\frac{dT}{d\tau}-1\right)=\csch^2{(nT-\tau)}\left(n\frac{dT}{d\tau}-1\right). 
\label{eq:deriv}
\ee
By using the identity $\csch^2{x}=\coth^2{x}-1$ and substituting~\eqref{eq:existB} into~\eqref{eq:deriv}, then using~\eqref{eq:kap2}, we obtain the result.

These superstable points on branches of solutions corresponding to different values of $n$ all occur at the same value of $T$. To see this note from~\eqref{eq:kap2} that $\tau-nT=\coth^{-1}{(\kappa/2)}$. Using the oddness of $\coth$ and substituting this into~\eqref{eq:existB} we get
\be
   \coth{[T-\coth^{-1}{(\kappa/2)}]}=\kappa/2
\ee
and, hence, $T=2\coth^{-1}{(\kappa/2)}$, which is independent of $n$. \qed
\end{enumerate}
In summary, if $\gamma$ varies monotonically along a branch of solutions (which we find to always occur), the solution is stable when 
$0<\gamma<1$ and superstable when $\gamma=1$. Moreover, it becomes unstable as $\gamma$ increases through $(n+1)/n$, where an eigenvalue increases through~1; generically, this is a saddle-node bifurcation. 

\begin{figure}
\begin{center}
\includegraphics[width=12.5cm]{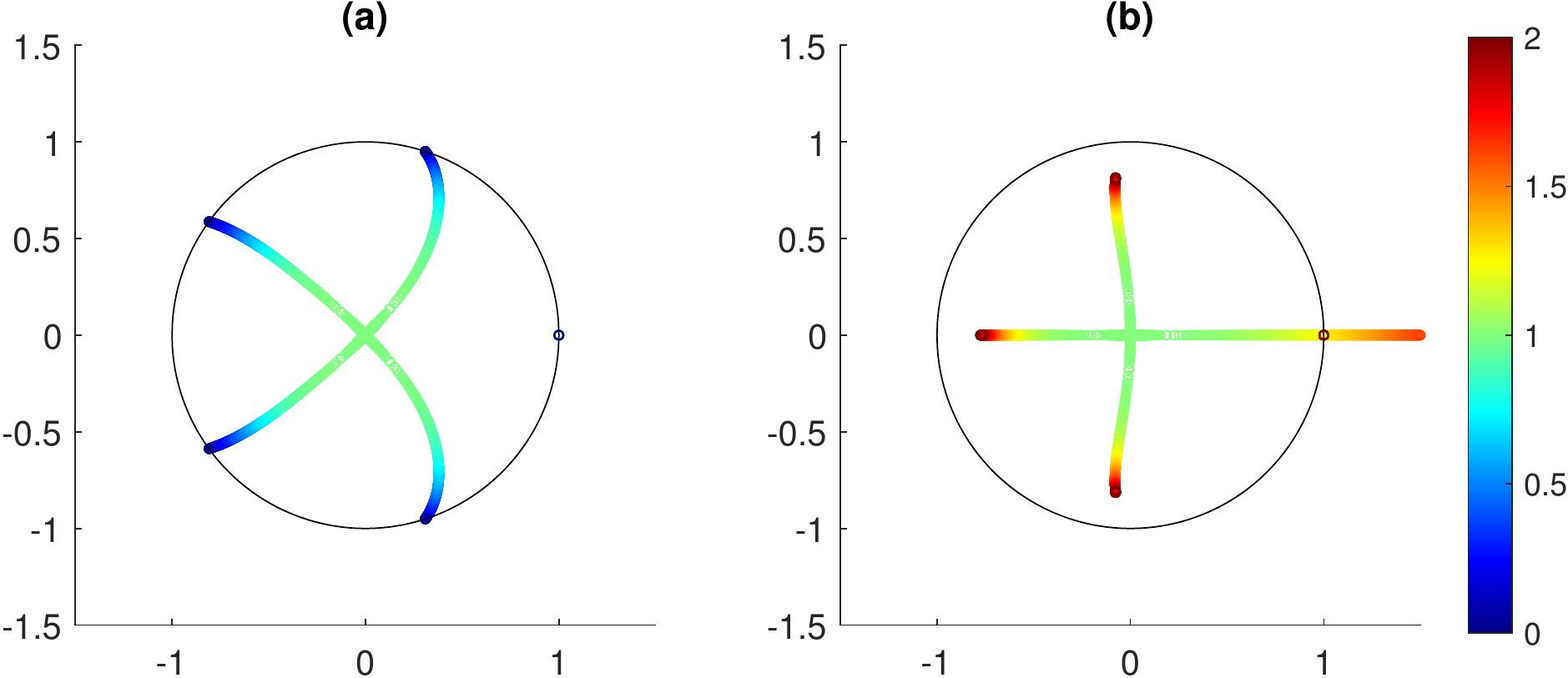}
\caption{Roots in the complex plane (with the unit circle shown) of the function $g(\lambda)$ from~\eqref{eq:defg} for $n=4$ as its parameter $\gamma$ (represented by the colour bar) is varied for
$0 \leq \gamma < 1$ (a) and for $1 < \gamma \leq 2$ (b).}
\label{fig:multip}
\end{center}
\end{figure}

Figure~\ref{fig:multip} illustrates Proposition~\ref{prop:g} by showing the typical behaviour of the roots of the function $g(\lambda)$ in the complex plane as the parameter $\gamma$ is varied; specifically, for the case $n=4$ of a periodic solution with five spikes in the delay interval $(-\tau,0]$. The value of $\gamma$ is represented by color in Fig.~\ref{fig:multip} and panel (a) shows the roots for $0<\gamma<1$. For $\gamma=0$ the roots of $g(\lambda)$ are the fifth roots of unity. As $\gamma$ increases, the four roots of the factor $h(\lambda)$ decrease in magnitude and their arguments tend to those of the fourth roots of $-1$, that is, to~$\pi/4,3\pi/4,5\pi/4$ and $7\pi/4$; at $\gamma = 1$ there is a quadruple root $0$ of $h$ and of $g$ and the periodic orbit is superstable. This behavior can be understood by considering the function $h(\lambda)$. Since $|\lambda|\to 0$ as $\gamma\to 1$, we can approximate $h(\lambda)$ by $\lambda^n+1-\gamma$, showing that its roots $\lambda$ tend to the $n$th roots of $\gamma-1$, which is negative here. Figure~\ref{fig:multip}(b) shows the behaviour of the roots of $g(\lambda)$ as $\gamma$ is increased from 1 in the interval $1 < \gamma \leq 2$. The four roots of $h(\lambda)$ leave the origin with arguments close to those of the fourth roots of $+1$, i.e.~$0,\pi/2,\pi$ and $3\pi/2$. This can also be understood with the argument above, since now $\gamma-1$ is positive. Moreover, we see that one real multiplier leaves the unit circle through $1$ as $\gamma$ increases through $(n+1)/n=5/4$. Note that the trivial root $1$ of $g(\lambda)$ remains unchanged throughout in Fig.~\ref{fig:multip}.

\subsection{Branches of periodic solutions for negative current}
\label{sec:numerA}

\begin{figure}
\begin{center}
\includegraphics[width=11.5cm]{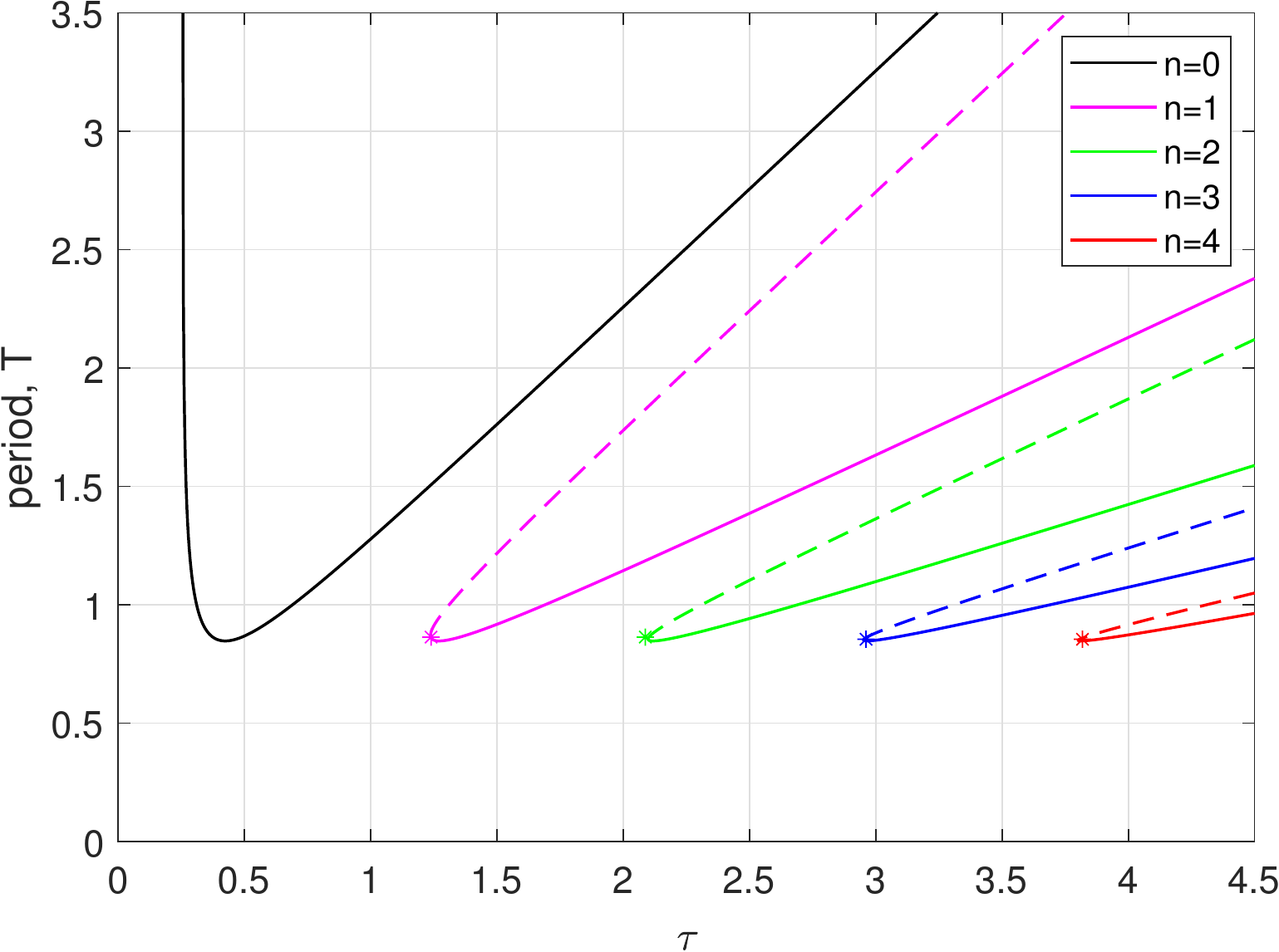}
\caption{Branches of periodic solutions for $n=0,1,2,3,4$ (see legend) of~\eqref{eq:dth} with negative current $I$, shown by their period $T$ as a function of the delay $\tau$. Periodic solutions are stable along solid curves, unstable along dashed curves, and superstable at the minima of $T$; stars indicate saddle-node bifurcations. Here, $I = -I_m^2 = -1$ and $\kappa=5$.}
\label{fig:pertauA}
\end{center}
\end{figure}

Figure~\ref{fig:pertauA} shows solution branches of Eq.~\eqref{eq:dth} with $\kappa=5$ (and $I_m = 1$) for different values of $n$ as indicated. Each branch represents a periodic solution with $n+1$ spikes in the delay interval, represented by its period $T$ as a function of $\tau$, with stability indicated as determined in the previous section. The single-spike solution for $n = 0$ is stable throughout. Branches for $n\geq1$ emerge in pairs, one stable and one unstable, at saddle-node bifurcations as $\tau$ is increased; the different stable solutions coexist, leading to an increasing level of multistability. Similar figures for other excitable systems with delayed feedback appear in~\cite{wedslo21,ruskra20,yanper09}. 

The solution branch for $n=0$ is special: it is always stable (when it exists) and the period approaches infinity as $\tau$ is decreased sufficiently. This corresponds to a homoclinic bifurcation to the unstable fixed point $\theta=2\tan^{-1}{(1)}$ at a finite value $\tau^*$.  To understand this, imagine that $\theta(0)=\pi$ on this homoclinic orbit. Under backwards time, $\theta(t)$ will approach $2\tan^{-1}{(1)}$ as $t\to -\infty$. In forwards time, at time $t=\tau^*$ the impulse from the delayed feedback will precisely make the neuron jump to $\theta=2\tan^{-1}{(1)}$, where it will then stay forever.  This condition is described by
\be
   1=\tan{\left(\frac{\theta(\tau^-)}{2}\right)}+\kappa
\ee
where $\theta(\tau^-)=-2\tan^{-1}{\left[\coth{(\tau)}\right]}$, 
that is, by
\be
   1=\kappa-\coth{(\tau)} 
\label{eq:hom}
\ee
For $\kappa=5$ we obtain the solution $\tau^*=\coth^{-1}{(\kappa-1)}\approx 0.25541$, which is the position of the vertical asymptote in Fig.~\ref{fig:pertauA}. The minimum in $T$ of this curve for $n=0$ occurs at $(\tau, T)=(\overline{T}/2,\overline{T}) \approx (0.42365, 0.8473)$, where $\overline{T} = 2 \coth^{-1}{(\kappa/2)}$ is the minimum period.

To determine the branches for $n \geq 1$ it is not necessary to solve~\eqref{eq:existB} for additional values of $n$. Instead, they can be found and plotted by using the reappearance of periodic solutions of DDEs with fixed delay~\cite{yanper09}, which implies that the branches in Fig.~\ref{fig:pertauA} are images of one another under a similarity transformation. Suppose that a DDE (with fixed delay) has a periodic solution with period $T_0$ for a time delay $\tau=\tau_0$. This same periodic solution is then also a solution for delay $\tau=\tau_0+mT_0$ where $m$ is an integer. Thus, a particular branch of periodic solutions of~\eqref{eq:dth} reappears at higher (and lower) values of the delay. 

We can express all of these branches parametrically. The branch for $n=0$, for which $\tau<T(\tau)$, is referred to as the {\em primary branch}. It is given explicitly as a function of $\tau$ by
\be
  T(\tau)=\tau+\coth^{-1}{(\kappa-\coth{\tau})}. 
\label{eq:prim}
\ee
The branches for $n \geq 1$ are referred to as {\em secondary branches}.
Letting $s$ be a parameter, where $\coth^{-1}{(\kappa-1)}<s<\infty$, the $n$th branch is of the form
\be
   (\tau,T)=(s+nT(s),T(s)). 
\label{eq:param}
\ee
The reappearance of periodic solutions can also be seen explicitly from~\eqref{eq:existB}. Suppose for some $n$ and $\tau=\tau_0$ there is a periodic solution with period $T_0$, i.e.
\[
   \coth{[(n+1)T_0-\tau_0]}=\kappa+\coth{(nT_0-\tau_0)}.
\]
Then it follows that
\[ \coth{[(n+m+1)T_0-(\tau_0+mT_0)]}=\kappa+\coth{[(n+m)T_0-(\tau_0+mT_0)]},
\]
where $m$ is any integer; hence, there is also a solution of period $T_0$ at delay $\tau=\tau_0+mT_0$ with $n+m$ past firing times in the interval $(-\tau,0)$.
This reappearance also explains why all superstable points occur at the same value $\overline{T}$ of $T$, as described in Sec.~\ref{sec:stab1}: they are the images of the minimum at $(\tau,T)= (\overline{T}/2,\overline{T})$ on the primary branch. It follows that the minimum on the $n=1$ branch is at $(\tau,T)=(3\overline{T}/2,\overline{T})$, on the $n=2$ branch at $(\tau,T)=(5\overline{T}/2,\overline{T})$, and on the $n$th branch at $(\tau,T)=((2n+1)\overline{T}/2,\overline{T})$. Thus, we conclude that the minimum on the $n$th branch occurs at its (nontransverse) intersection point with the line $T=2\tau/(2n+1)$. 

\begin{figure}
\begin{center}
\includegraphics[width=7cm]{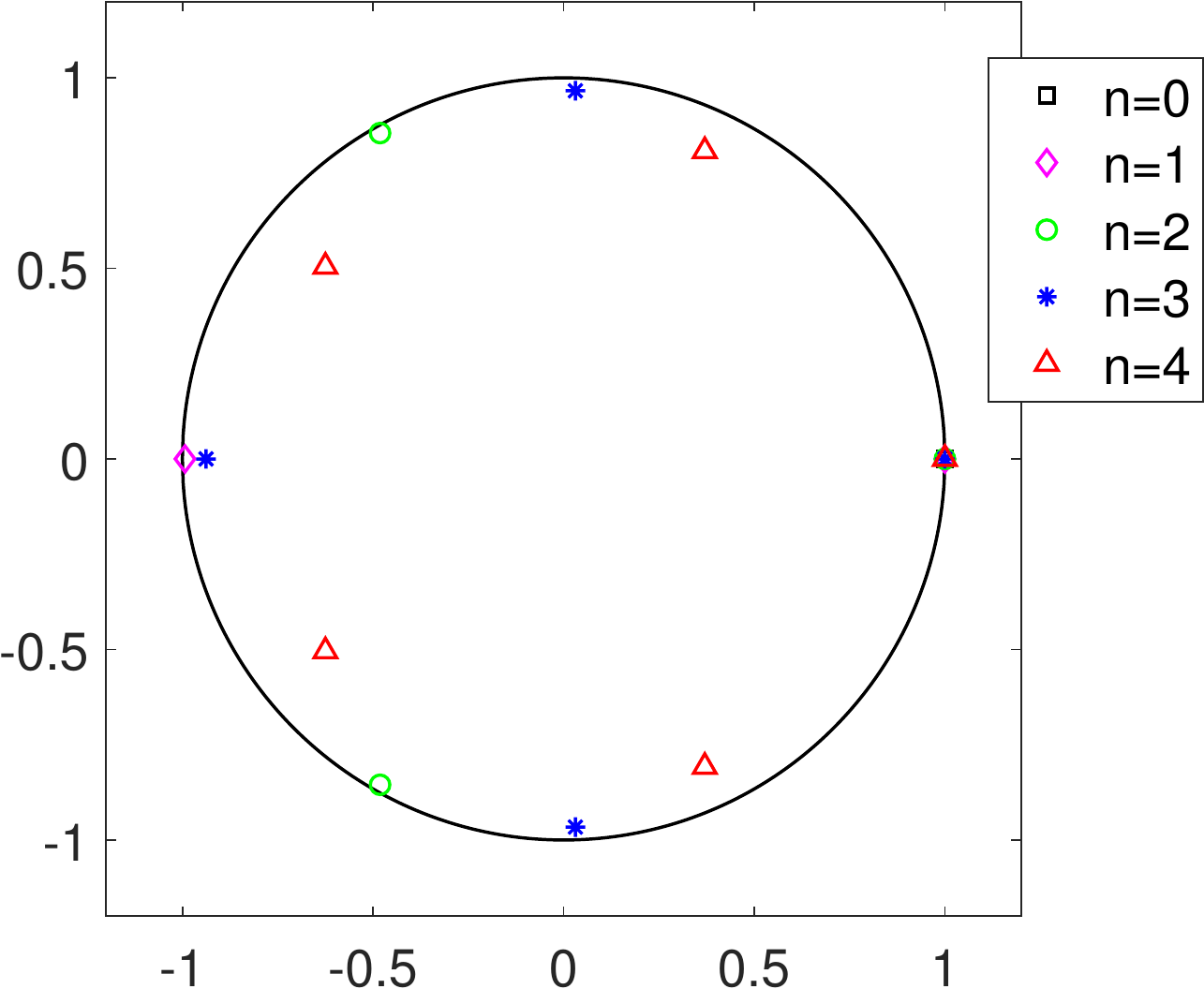}
\caption{The Floquet multipliers of the simultaneously stable periodic solutions with $n=0,1,2,3,4$ (see legend) for $\tau=4$ in Fig.~\ref{fig:pertauA}, as computed from~\eqref{eq:defg}. 
Here, $I = -I_m^2 = -1$ and $\kappa=5$. 
Compare with Fig.~\ref{fig:multip}.}
\label{fig:roots}
\end{center}
\end{figure}

We conclude this section by considering the stability of simultaneously stable periodic solutions. Figure~\ref{fig:roots} shows the sets of Floquet multipliers for the different stable solutions given by~\eqref{eq:existB} with $n=0,1,2,3,4$ at $\tau=4$ in Fig.~\ref{fig:pertauA}. As discussed in Sec.~\ref{sec:stab1}, these multipliers are roots of the polynomial $g(\lambda)$ from \eqref{eq:defg} and can, thus, easily be found numerically for any value of $\gamma$ (see~\eqref{eq:gamma}) which in turn depends on $\kappa$ and $\tau$. For solutions such as those in Fig.~\ref{fig:pertauA}, the parameter $\gamma$ tends to zero from above as one moves along a stable branch away from the saddle-node bifurcation where it is created. According to Proposition~\ref{prop:g}(i), the respective Floquet multipliers approach the unit circle from within as $\gamma\to 0$, namely at roots of unity; see also Fig.~\ref{fig:multip}. This explains why the multipliers for smaller $n$ are closer to the unit circle, which is certainly the case for $n=0,1,2,3$ in Fig.~\ref{fig:roots}, meaning that the stability of the corresponding periodic solutions is already quite weak. This analytical result for the theta neuron~\eqref{eq:dth} with Dirac delta function is in agreement with, and may serve as an explanation for, the observation in other contexts~\cite{RomeiraNSR16,TerrienPRE21} of only weakly stable solutions with $n$ spikes or pulses in the delay interval, whose Floquet multipliers are near roots of unity.

\subsection{Bifurcation curves in the $(\tau,\kappa)$-plane for negative current}
\label{sec:bifdiagExc}

\begin{figure}
\begin{center}
\includegraphics[width=9.5cm]{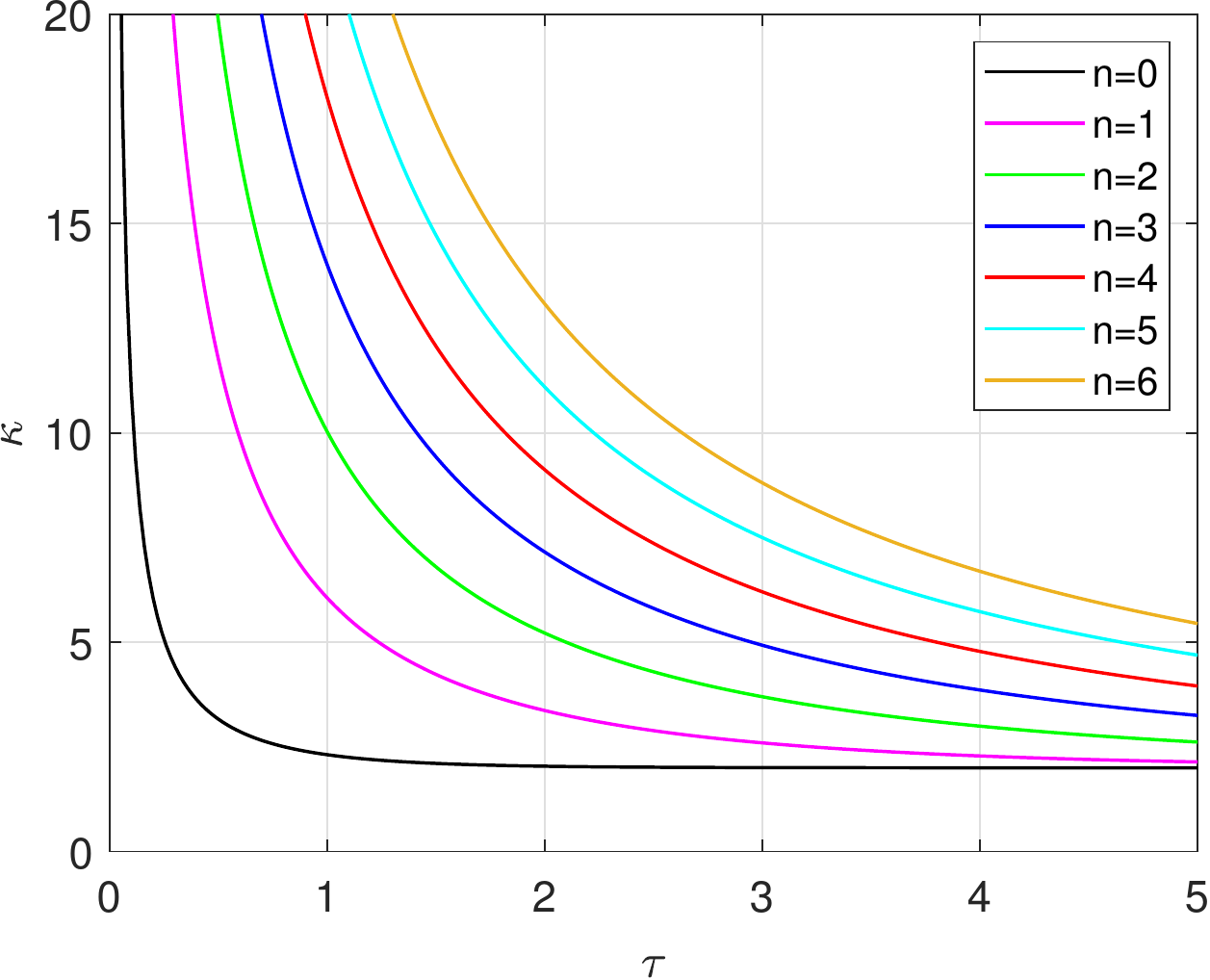}
\caption{Bifurcation curves in the $(\tau,\kappa)$-plane for $n=0,1,2,3,4,5,6$ (see legend) of~\eqref{eq:dth} with negative current $I$; compare with Fig.~\ref{fig:pertauA}. The curve for $n=0$ (black) is the locus of homoclinic bifurcation of the primary branch; no periodic solutions exist below and to the left of this curve. The curves for $n=1,2,3,4,5,6$ are the saddle-node bifurcations of the secondary branches. Here, $I = -I_m^2 = -1$.}
\label{fig:kaptauA}
\end{center}
\end{figure}

The bifurcation of the branches in Fig.~\ref{fig:pertauA} can be continued in the additional parameter $\kappa$. Figure~\ref{fig:kaptauA} shows in the $(\tau,\kappa)$-plane the curves of homoclinic bifurcation on the primary branch, given by~\eqref{eq:hom}, and of saddle-node bifurcations on the secondary branches for $n=1,2,3,4,5,6$. Generally, such loci of bifurcations need to be continued numerically with standard numerical algorithms~\cite{kuz04,auto,NewDDEBiftool}; see also Sec.~\ref{sec:smooth}. However, for the theta neuron with delta feedback \eqref{eq:dth} they can actually be found analytically. Writing the value of $\tau$ on the $n$th branch as $\tau_n=\tau_0+nT(\tau_0)$ where $\tau_0$ is the value of $\tau$ on the primary branch, a saddle-node bifurcation occurs when $d\tau_n/d\tau_0=0$~\cite{ruskra20}, i.e. when~$0=1+nT'(\tau_0)$, or $T'(\tau_0)=-1/n$. Differentiating~\eqref{eq:prim} we have
\[
   T'(\tau_0)=1+\frac{\coth^2{\tau_0}-1}{1-(\kappa-\coth{\tau_0})^2}.
\]
Setting this $T'(\tau_0) = -1/n$ and then solving the resulting quadratic equation for $\coth{\tau_0}$ we obtain
\be
   \coth{\left(\tau_0^{(n)}\right)}=\kappa(1+n)-\sqrt{1+\kappa^2(n^2+n)}, \label{eq:cothtau}
\ee
where we have taken the physically meaningful of the two possible solutions.  
Substituting~\eqref{eq:cothtau} into~\eqref{eq:prim} we obtain the value of $T$ at which the saddle-node bifurcation on the $n$th branch occurs as
\be T^{(n)}=\coth^{-1}{\left[\kappa(1+n)-\sqrt{1+\kappa^2(n^2+n)}\right]}+\coth^{-1}{\left[\sqrt{1+\kappa^2(n^2+n)}-\kappa n\right]}.
\label{eq:valueSN_T}
\ee

\setlength{\unitlength}{1mm}

\begin{figure}
\begin{center}
\begin{picture}(80,65)(0,0)
\put(0,0){\includegraphics[width=8cm]{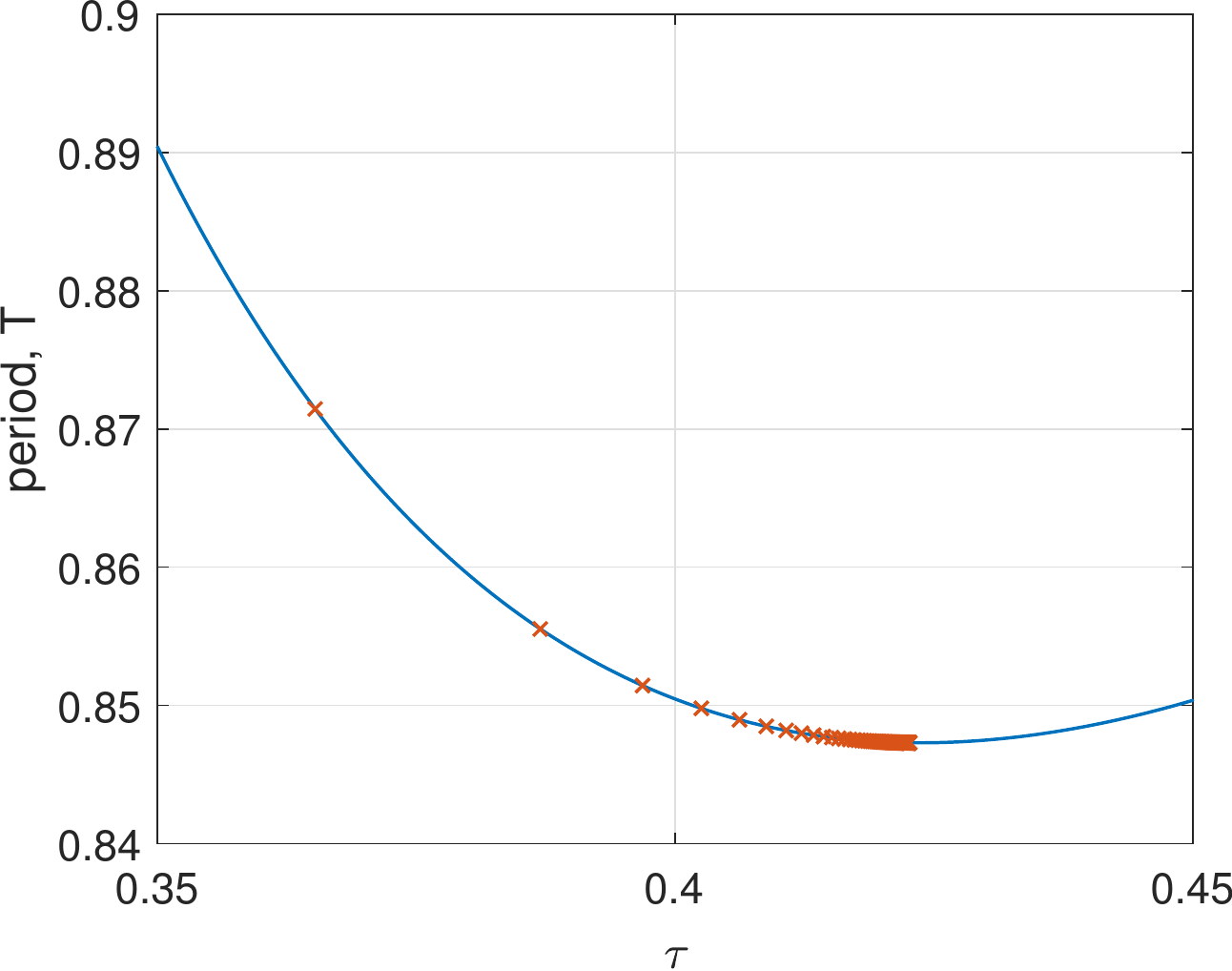}}
\put(21,37){$(\tau_0^{(1)},T^{(1)})$}
\put(35.7,23){$(\tau_0^{(2)},T^{(2)})$}
\put(59,10.5){$(\overline{T}/2,\overline{T})$}
\put(58.48,14){$\circ$}
\end{picture}
\end{center}
\caption{The points $(\tau_0^{(n)},T^{(n)})$ (red crosses) on the primary branch of periodic solutions (blue curve) converge to its minimum at $(\tau,T) = (\overline{T}/2,\overline{T})$; these points map to the points of saddle-node bifurcation on the $n$th branch at $(\tau_0^{(n)}+nT^{(n)},T^{(n)})$; compare with Fig.~\ref{fig:pertauA}. Here, $I = -I_m^2 = -1$ and $\kappa=5$.}
\label{fig:snmap}
\end{figure}
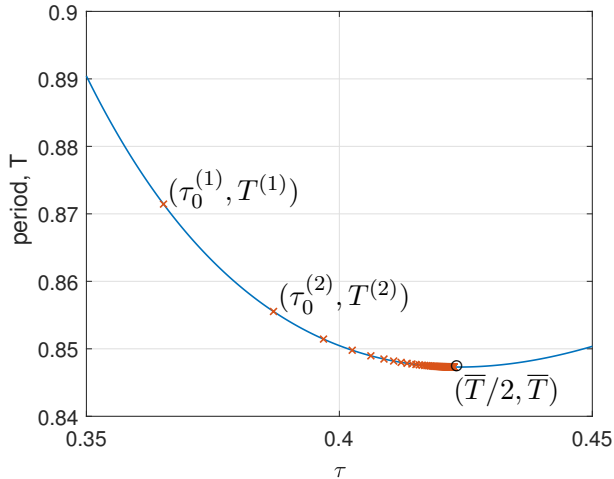

In Fig.~\ref{fig:snmap} the first 100 of the points $T^{(n)}$ for $\kappa=5$ in \eqref{eq:valueSN_T} are plotted on the primary branch. This illustrates that $T^{(n)}$ coverges to the value $\overline{T}$ at the minimum of the primary branch. Indeed,
\[
   \lim_{n\to\infty}\kappa(1+n)-\sqrt{1+\kappa^2(n^2+n)}=\kappa/2
\]
and
\[
   \lim_{n\to\infty}T^{(n)}=2\coth^{-1}{(\kappa/2)} = \overline{T}.
\]
The value of $\tau$ at which the saddle-node bifurcation of the $n$th branch occurs is $\tau^{(n)}=\tau_0^{(n)}+nT^{(n)}$; this expression depends on $\kappa$ and the curves in Fig.~\ref{fig:kaptauA}  were obtained by plotting $\tau^{(n)}$ as a function of $\kappa$ for various $n$s indicated.

\section{Oscillating theta neuron for positive current}
\label{sec:pos}

We now consider~\eqref{eq:dth} with a positive current $I=I_p^2>0$, which means that the uncoupled
theta neuron is producing regular spikes. We first consider the existence and stability of periodic orbits in the presence of delayed self-coupling. Subsequently, we consider the associated branches of periodic solutions and their bifurcations for the two sub-cases of excitatory self-coupling for $\kappa>0$, and of inhibitory self-coupling for $\kappa<0$ --- showing that they are related by an explicit transformation.

\subsection{Existence of periodic solutions for positive current}
\label{sec:posex}

Our starting point here is again the solution of~\eqref{eq:QIF} with $\kappa=0$, which for with $I=I_p^2>0$ is
\[
   V(t)=I_p\tan{[I_pt+\tan^{-1}{(V(0)/I_p)}]}.
\]
With the transformation $V=\tan{(\theta/2)}$ we have
\be \theta(t)=2\tan^{-1}\left[I_p\tan{\left(I_pt+\tan^{-1}\left(\frac{\tan{(\frac{\theta(0)}{2})}}{I_p}\right)\right)}\right]. 
\label{eq:solI}
\ee
To find periodic solutions for $\kappa\neq0$ we assume, as before, that the neuron fires at $t=0$ and there are $n$ past firing times in $(-\tau,0)$. The delta function acts at time $\tau-nT$, at which point
\[ \theta(\tau-nT^-)=2\tan^{-1}\left[I_p\tan{\left(I_p(\tau-nT)+\frac{\pi}{2}\right)}\right].
\]
The delta function acts to move the phase to $\theta(\tau-nT^+)$ where
\[ \tan{\left(\frac{\theta(\tau-nT^+)}{2}\right)}=\tan{\left(\frac{\theta(\tau-nT^-)}{2}\right)}+\kappa=I_p\tan{\left[I_p(\tau-nT)+\frac{\pi}{2}\right]}+\kappa.
\]
The neuron fires (i.e.~$\theta$ reaches $\pi$) after another time $\Delta$ 
where (using~\eqref{eq:solI})
\[ \pi=2\tan^{-1}\left[I_p\tan{\left(I_p\Delta+\tan^{-1}\left(\frac{\tan{(\frac{\theta(\tau-nT^+)}{2})}}{I_p}\right)\right)}\right],
\]
i.e.~when 
\[ I_p\Delta+\tan^{-1}\left[\frac{\tan{(\frac{\theta(\tau-nT^+)}{2})}}{I_p}\right]=\frac{\pi}{2},
\]
which gives 
\[ \Delta=\frac{1}{I_p}\left[\frac{\pi}{2}-\tan^{-1}{\left[\kappa/I_p+\tan{\left(I_p(\tau-nT)+\frac{\pi}{2}\right)}\right]}\right].
\]
Since $\tau-nT+\Delta=T$ we have
\be (n+1)T=\tau+\frac{1}{I_p}\left[\frac{\pi}{2}-\tan^{-1}{\left[\kappa/I_p+\tan{\left(I_p(\tau-nT)+\frac{\pi}{2}\right)}\right]}\right].
\label{eq:exA}
\ee
As before, we can rescale either $I_p,\kappa$ or $\tau$ to be 1, so set $I_p=1$ and rewrite~\eqref{eq:exA} as
\be (n+1)T=\tau+\frac{\pi}{2}-\tan^{-1}{\left[\kappa+\tan{\left(\tau-nT+\frac{\pi}{2}\right)}\right]}. 
\label{eq:existA}
\ee
This is the equation relating the period $T$ to the other parameters $\tau$ and $\kappa$, for a given integer $n$.
As above, for $n=0$ expression~\eqref{eq:existA} gives $T$ explicitly as
\be T(\tau)=\tau+\frac{\pi}{2}-\tan^{-1}{\left[\kappa+\tan{\left(\tau+\frac{\pi}{2}\right)}\right]} = \tau+\frac{\pi}{2}-\tan^{-1}{\left[\kappa-\cot{\tau}\right]},
\label{eq:primA}
\ee
which is valid for $0\leq \tau\leq \pi$. Note that, when $I=I_p=1$, $d\theta/dt=2$ except at the times at which the feedback acts, which
simplifies the derivations below.

\subsection{Stability of periodic solutions for positive current}

Determining the stability for $I=I_p^2>0$ is also similar to the excitable case.
Assume again that the neuron has just fired at time $t_0$ and there are $n$ past firing times in $(-\tau,0)$. We wait $\tau-(t_0-t_{-n})$ until the delta function acts, which maps $\theta$ from
\[
  \theta(\tau-t_0+t_{-n}^-)=\pi+2(\tau-t_0+t_{-n})
\]
to
\[ \theta(\tau-t_0+t_{-n}^+)=2\tan^{-1}{\left[\tan{\left(\frac{\theta(\tau-t_0+t_{-n}^-)}{2}\right)}+\kappa\right]}.
\]
We then wait a time $\Delta$ until the neuron fires at time $t_1$, where
\[
   \pi=2\Delta+\theta(\tau-t_0+t_{-n}^+).
\]
Thus 
\[ t_1=t_0+(\tau-t_0+t_{-n})+\Delta=\tau+t_{-n}+\frac{\pi}{2}-\tan^{-1}{\left[\tan{\left(\tau-t_0+t_{-n}+\frac{\pi}{2}\right)}+\kappa\right]}
\]
and, in general,
\[ t_i=\tau+t_{i-n-1}+\frac{\pi}{2}-\tan^{-1}{\left[\tan{\left(\tau-t_{i-1}+t_{i-n-1}+\frac{\pi}{2}\right)}+\kappa\right]}.
\]
This is a map giving the next firing time, $t_i$, in terms of the previous ones, back to $t_{i-n-1}$. Perturbing the firing times defined by this equation as in Sec.~\ref{sec:stab1}, we obtain the same matrix $\J$ as in~\eqref{eq:J}. However, now for a periodic orbit with period $T$ we have
\be 
\gamma=\frac{\sec^2{(\tau-nT+\frac{\pi}{2})}}{1+\left[\tan{(\tau-nT+\frac{\pi}{2})}+\kappa\right]^2}=\frac{\csc^2{(\tau-nT)}}{1+\left[\kappa-\cot{(\tau-nT)}\right]^2},
\label{eq:stabA}
\ee
which is clearly positive for any value of $\kappa$.

\subsection{Branches of periodic orbits for positive current}

As in Sec.~\ref{sec:numerA}, the primary branch of periodic solutions is given by~\eqref{eq:primA} and the secondary branches are given parametrically by~\eqref{eq:param}. To plot and discuss these branches we need to distinguish the sub-cases $\kappa>0$ of excitatory and $\kappa<0$ of inhibitory self-coupling.

\subsubsection{The case of excitatory delayed self-coupling}

\begin{figure}
\begin{center}
\includegraphics[width=11.5cm]{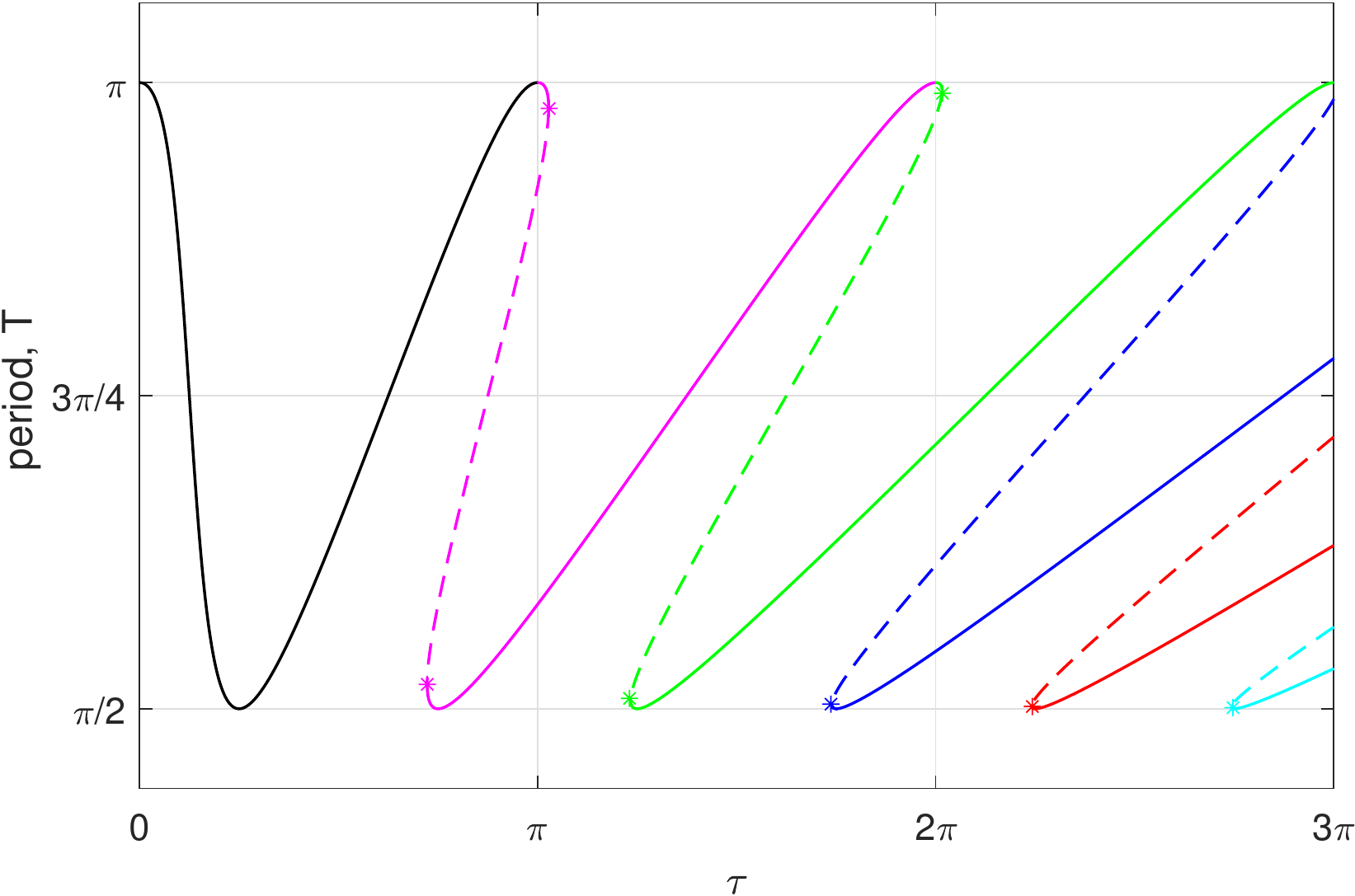}
\caption{Branches of periodic solutions for $n=0,1,2,3,4,5$ of~\eqref{eq:dth} with positive current $I$ and positive $\kappa$, showing their period $T$ as a function of the delay $\tau$. Periodic solutions are stable along solid curves, unstable along dashed curves, and superstable at the minima of $T$ and at the maxima of $T$ with $\tau > 0$; stars indicate saddle-node bifurcations. Here, $I = I_p^2 = 1$ and $\kappa=2$.}
\label{fig:pertauposA}
\end{center}
\end{figure}

Figure~\ref{fig:pertauposA} shows the branches of periodic solutions for $n=0,1,2,3,4,5$ with $I = I_p^2 = 1$ and $\kappa=2$ so that the self-coupling is excitatory. Notice that all branches now connect to form a single curve of periodic solutions; similar plots appear in~\cite{hasval12,yanper09}. However, different parts of this curve still correspond to different values of $n$, as is indicated by colour. The primary branch for $n=0$ is again given by a function $T(\tau)$ and it has a finite period throughout for positive current $I$; that is, there is no longer a homoclinic bifurcation. More specifically, the maximum period now occurs at
$\tau=0$ and is equal to that of an uncoupled oscillator; this is the case since
at firing, $\theta=\pi$ and in~\eqref{eq:dth} the term $1+\cos{\theta}$ is equal to zero, so the feedback can have no effect. The period of this free oscillation is $\pi/I_p=\pi$. For the chosen value of $\kappa$ the primary branch is entirely stable. On the other hand, all secondary branches have two saddle-node bifurcations on them, with an unstable middle sub-branch. Notice also that the secondary branches are increasingly tilted to the right as $n$ increases, leading again to an increasing level of multistability with increasing $\tau$.  

A periodic orbit is superstable when $dT/d\tau=0$, and this again occurs at the minima of the period but now also at the maximum period (with $\tau > 0$) where neighbouring branches meet. Differentiating~\eqref{eq:primA} with respect to $\tau$, setting this to zero, and substituting the expression for $\tau$ back into~\eqref{eq:primA}, we find that the minimum period on the primary branch  occurs at $(\tau, T) = (\overline{T}/2,\overline{T})$ where now $\overline{T} = 2\cot^{-1}(\kappa/2)$;  the minima of the secondary branches are therefore at $(\tau, T) = ((n+1/2)\overline{T}, \overline{T})$.
(For this value of $\kappa$ we have $\overline{T} =\pi/2$.)
The primary branch exists for $0\leq\tau\leq \pi$, so the transitions between branches occur when $\tau$ is a multiple of $\pi$; these are also the locations of the maxima where the periodic orbit is also superstable.

\subsubsection{The case of inhibitory delayed self-coupling}

\begin{figure}
\begin{center}
\includegraphics[width=11.5cm]{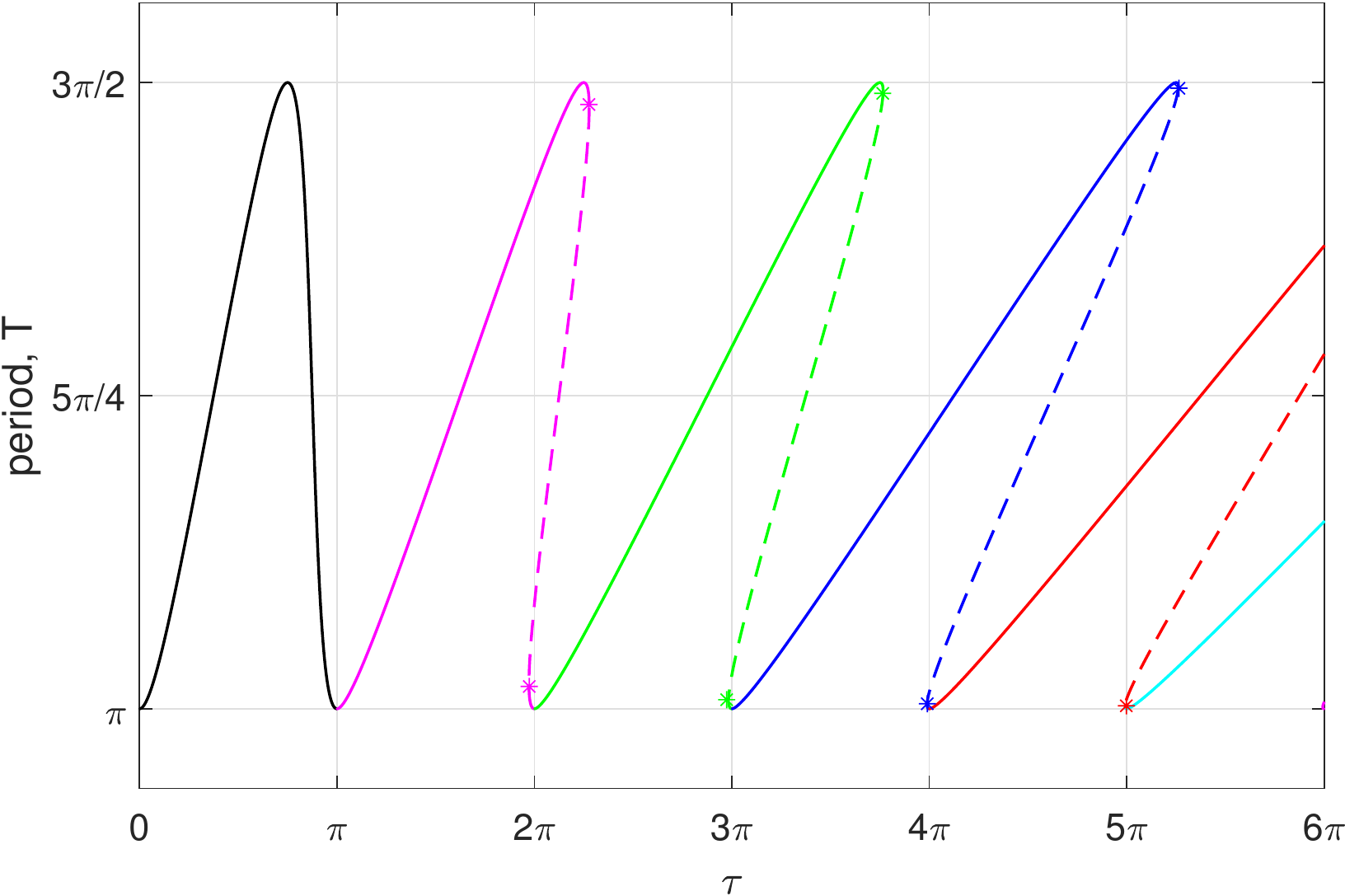}
\caption{Branches of periodic solutions for $n=0,1,2,3,4,5$ of~\eqref{eq:dth} with positive current $I$ and negative $\kappa$, showing their period $T$ as a function of the delay $\tau$. Periodic solutions are stable along solid curves, unstable along dashed curves, and superstable at the minima of $T$; stars indicate saddle-node bifurcations. Here, $I = I_p^2 = 1$ and $\kappa=-2$; compare with Fig.~\ref{fig:pertauposA}.}
\label{fig:pertaunegA}
\end{center}
\end{figure}

Figure~\ref{fig:pertaunegA} shows the shows the branches of periodic solutions for $n=0,1,2,3,4,5$ with $I = I_p^2 = 1$ and an inhibitory self-coupling with $\kappa=-2$. At $\tau=0$ the period is again $\pi$ but this is now the minimum period possible; this makes sense since inhibitory feedback can only increase the period. The primary branch is again entirely stable, while the secondary branches each have two points of saddle-node bifurcation that delimit a sub-branch where the periodic orbit is unstable.

The branches in Fig.~\ref{fig:pertaunegA} have been plotted by evaluating~\eqref{eq:primA} and~\eqref{eq:param} for $\kappa < 0$. However, they can also be obtained directly from those in  Figs.~\ref{fig:pertauposA} via the following interesting geometric relationship between the branches of periodic solutions for $\kappa > 0$ and for $\kappa < 0$.

\begin{prop}[Transformation between excitatory and inhibitory cases]
\verb+ + \\[-3mm]  
\label{prop:excinh}

\noindent
For positive current $I=1$ and given $n \geq 0$, the $n$th branch of periodic orbits of~\eqref{eq:existA} for $\kappa = K > 0$ maps to the $n$th branch for $\kappa =  -K <  0$ under rotation over $\pi$ about the point
\be
   (\tau_c,T_c)=\left((n+1/2)\pi,\pi\right),
\label{eq:centerofrotationA}
\ee
and vice versa; see Figs. ~\ref{fig:pertauposA} and~\ref{fig:pertaunegA}.

\end{prop}

\noindent
{\bf Proof:}
Suppose that $(\tau,T)$ is a point on the $n$th branch for $\kappa = K > 0$, 
i.e.~$(\tau,T)$ satisfies
\be (n+1)T=\tau+\frac{\pi}{2}-\tan^{-1}{\left[K+\tan{\left(\tau-nT+\frac{\pi}{2}\right)}\right]}.
\label{eq:rot}
\ee 
Rotating this point through $\pi$ radians about the center of rotation $(\tau_c,T_c)$ from \eqref{eq:centerofrotationA} 
gives the new point
\[
   (\widehat{\tau},\widehat{T}) = (2(n+1/2)\pi-\tau,2\pi-T).
\]
Solving for $\tau$ and $T$ and substituting back into~\eqref{eq:rot}, we see that
$(\widehat{\tau},\widehat{T})$ satisfies
\[
(n+1)\widehat{T}=\widehat{\tau}+\frac{\pi}{2}-\tan^{-1} {\left[-K+\tan{\left(\widehat{\tau}-n\widehat{T}+\frac{\pi}{2}\right)}\right]}, 
\]
which is exactly of the form~\eqref{eq:rot} but now for $\kappa = -K$. \qed

\medskip

It follows immediately from Proposition~\ref{prop:excinh} that the minima in Fig.~\ref{fig:pertaunegA} all occur at $T=\pi$ when $\tau$ is a multiple of
$\pi$; similarly, the maximum on the $n$th branch is at $(\tau, T) = ((n+1/2)(2\pi+\overline{T})/2, 2\pi+\overline{T})$; recall here that $\overline{T} = 2 \cot^{-1}(\kappa/2)$ and that $\kappa$ is now negative.

\subsection{Bifurcation curves in the $(\tau,\kappa)$-plane for positive current}

To find the saddle-node bifurcations on the $n$th branch for $I = I_p^2 > 0$ we follow the analysis in Sec.~\ref{sec:numerA} of letting $\tau_0$ be the value of $\tau$ on the primary branch and solving $T'(\tau_0)=-1/n$, where $T(\tau)$ is given by~\eqref{eq:primA}. This now gives (generically) two values of $\tau_0$ for each $n$, given by
\be \tau_{0,\pm}^{(n)}=\cot^{-1}{\left[\kappa(n+1)\pm\sqrt{\kappa^2(n^2+n)-1}\right]}.
\label{eq:tausnpm}
\ee
Substituting these into~\eqref{eq:primA} gives the corresponding $T$-values
\be T_\pm^{(n)}=\cot^{-1}{\left[\kappa(n+1)\pm\sqrt{\kappa^2(n^2+n)-1}\right]}+\frac{\pi}{2}-\tan^{-1}{\left[\mp\sqrt{\kappa^2(n^2+n)-1}-\kappa n\right]}. \label{eq:Tsnpm}
\ee
We now consider first the case that $\kappa$ is positive. Then the upper sign corresponds to the saddle-node bifurcation with the larger period and the lower sign to that with the smaller period; see Fig.~\ref{fig:pertauposA}. The saddle-node bifurcations with $T_\pm^{(n)}$ occur at
\be
   \tau=\tau_{0,\pm}^{(n)}+nT_\pm^{(n)}. 
\label{eq:tausn}
\ee
Note that for the values of $\tau_{0,\pm}^{(n)}$ to be real we must 
have $\kappa^2(n^2+n)>1$.
So for a fixed $\kappa$ only branches with 
\be
   \frac{\sqrt{1+4/\kappa^2}-1}{2}<n
\label{eq:condcusp}
\ee
have saddle-node bifurcations on them. For $\kappa = 2$, as in Fig.~\ref{fig:pertauposA}, this is satisfied for all $n \geq 1$, which means that all secondary branches have a pair of saddle-node bifurcations on them. (The same is true for Fig.~\ref{fig:pertaunegA}.)

Equation~\eqref{eq:tausn} (via~\eqref{eq:Tsnpm}) provides, for each choice of the sign, the value of $\tau$ at the respective saddle-node bifurcation of the $n$th branch as an explicit function of $\kappa$. For each $n$ satisfying \eqref{eq:condcusp} there are two curves in the $(\tau,\kappa)$-plane and they coincide at cusp points when $\kappa^2(n^2+n)=1$, that is, at 
\be (\tau,\kappa)=\left[(n+1)\cot^{-1}\left(\sqrt{\frac{n+1}{n}}\right)+\frac{n\pi}{2}+n\tan^{-1}\left(\sqrt{\frac{n}{n+1}}\right),\frac{1}{\sqrt{n^2+n}}\right].
\ee

\begin{figure}
\begin{center}
\includegraphics[width=10cm]{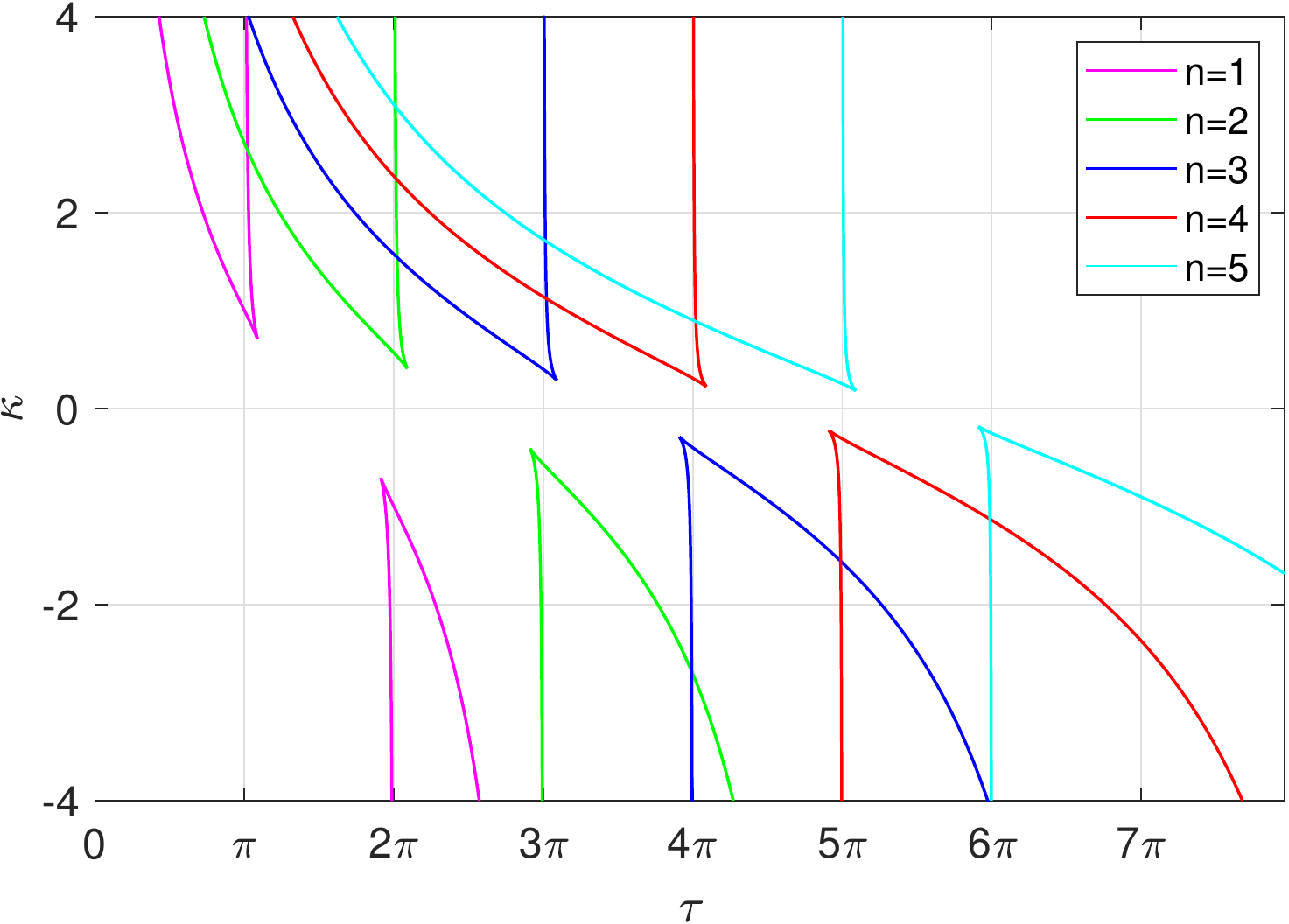}
\caption{Curves of saddle-node bifurcations of periodic orbits in the $(\tau,\kappa)$-plane for $n=1,2,3,4,5$ (see legend) of~\eqref{eq:dth} with positive current $I=I_p^2 =1$. The curves with $\kappa>0$ are for excitatory self-feedback, and those with $\kappa<0$ are for inhibitory self-feedback. Compare with Figs.~\ref{fig:pertauposA} and~\ref{fig:pertaunegA}, respectively.}
\label{fig:kaptaupos}
\end{center}
\end{figure}

Figure~\ref{fig:kaptaupos} shows the curves of saddle-node bifurcations for $n=1,2,3,4,5$ of periodic orbits of~\eqref{eq:existA} with $I = I_p^2 = 1$. 
The case of excitatory self-feedback is shown in the upper half of the $(\tau,\kappa)$-plane where $\kappa>0$. Notice that the respective cusp points, where the two curves defined by~\eqref{eq:tausn} come together, are minima here.

Figure~\ref{fig:kaptaupos} also shows the curves of saddle-node bifurcations for $n=1,2,3,4,5$ for the case of inhibitory coupling, namely in the lower half of the $(\tau,\kappa)$-plane where $\kappa<0$. These curves can be obtained from Eq.~\eqref{eq:tausn} (via~\eqref{eq:Tsnpm}), now for negative $\kappa$. However, Proposition~\ref{prop:excinh} implies a transformation also of the loci of saddle-node bifurcations on the $n$th branch of periodic solutions when the sign of $\kappa$ is changed. Specifically, these two loci are each other's images under rotation by $\pi$ about the point
\be
   (\tau_c,\kappa_c)=\left((n + 1/2)\pi,0\right).
\label{eq:centerofrotation}
\ee
It follows that the cusp points for negative $\kappa$ are at 
\be (\tau,\kappa)=\left[\left(1+\frac{3n}{2}\right)\pi-(n+1)\cot^{-1}\left(\sqrt{\frac{n+1}{n}}\right)-n\tan^{-1}\left(\sqrt{\frac{n}{n+1}}\right),\frac{-1}{\sqrt{n^2+n}}\right].
\ee

\section{Theta neuron with smooth self-feedback}
\label{sec:smooth}

We now consider the model of a theta neuron with smooth feedback given by 
\be
   \frac{d\theta}{dt}=1-\cos{\theta}+(1+\cos{\theta})\left[I+\kappa P(\theta(t-\tau))\right] 
\label{eq:dthsm}
\ee
where $\tau$ is the delay and
\be
   P(\theta)=(8/63)(1-\cos{\theta})^5 \label{eq:ptheta}
\ee
is a pulsatile function centred at $\theta=\pi$; here the factor $8/63$ is a normalisation ensuring that $\int_0^{2\pi} P(\theta)d\theta=2\pi$. 

\begin{figure}
\begin{center}
\includegraphics[width=8cm]{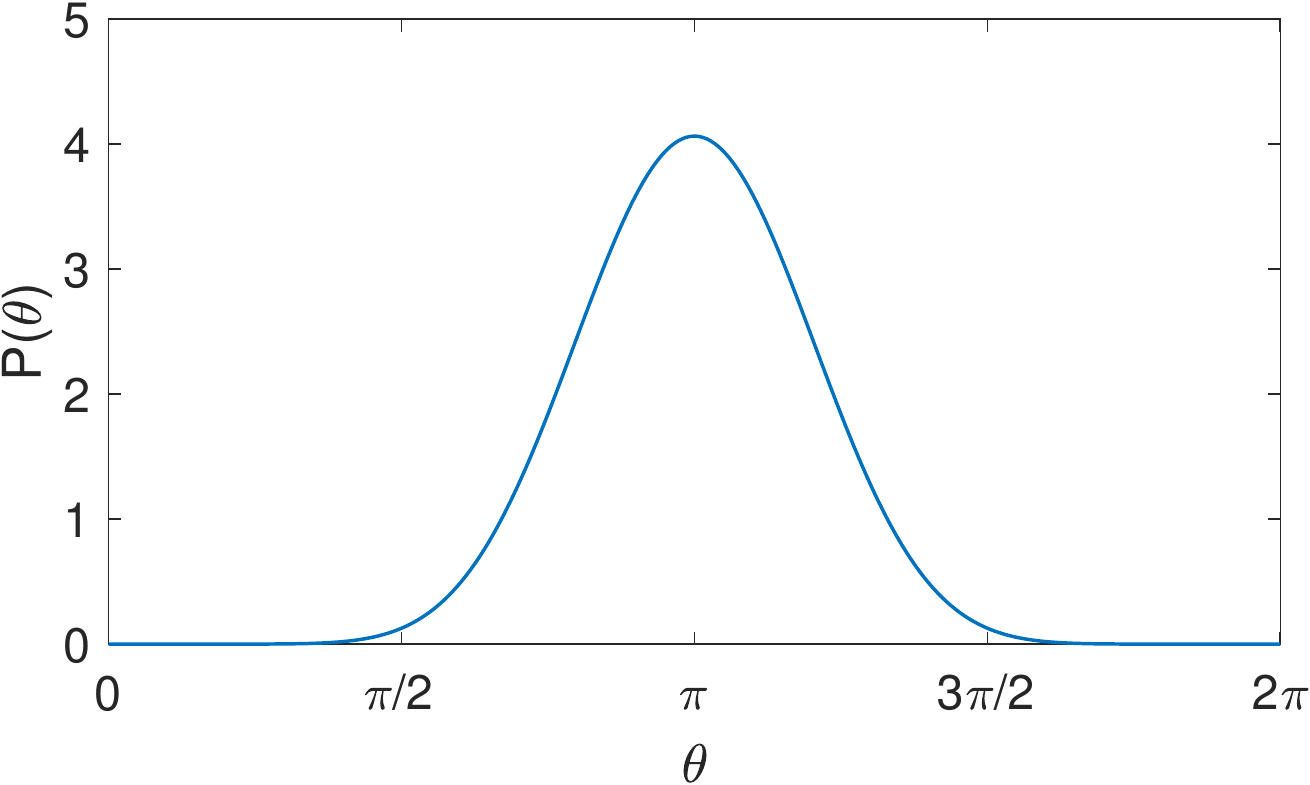}
\caption{The pulsatile function $P(\theta)$ given in~\eqref{eq:ptheta} as used 
in the theta neuron with smooth self-feedback~\eqref{eq:dthsm}.}
\label{fig:poftheta}
\end{center}
\end{figure}

The function $P(\theta)$ is shown in Fig.~\ref{fig:poftheta}. Its pulse is quite broad, certainly compared to the Dirac delta function that models the instantaneous delay in system \eqref{eq:dth}. Such forms of smooth coupling have been considered elsewhere~\cite{lai14A,lukbar13,chahat17}, although in infinite networks and without delays. In contrast, the mathematically similar Kuramoto model of coupled phase ocillators with delayed coupling has been well-studied~\cite{earstr03,yeustr99,leeott09}. 

We study here the theta neuron with smooth feedback~\eqref{eq:dthsm}--\eqref{eq:ptheta}, or \emph{smooth theta neuron} for short, to investigate the validity of the results we found for the case of instantaneous feedback in system \eqref{eq:dth}. A significant difference between these two systems is that $P$ is a function of $\theta$, not of time. So if $\theta$ increases through $\pi$ with a speed bounded away from zero, the function $P(\theta(t))$ will be pulse-like in time with a maximum at the time at which $\theta=\pi$. However, if $\theta$ {\em decreases} through $\pi$ (as a result of inhibitory coupling, for example) the neuron will then emit a spurious pulse. 

The smooth theta neuron~\eqref{eq:dthsm}--\eqref{eq:ptheta} is a delay differential equation (DDE) with a single fixed delay. As such, it is an infinite-dimensional dynamical system whose equilibria and periodic orbits must be expected to undergo  standard bifurcations; see, for example, \cite{Diekmann1995,Guo2013,HaleVerduyn1993,HalSmith2011}. Periodic solutions of DDEs such as~\eqref{eq:dthsm}--\eqref{eq:ptheta} are not known analytically but must be found with numerical methods \cite{KraSie2022,ROO07}. To find (stable) periodic solutions we numerically integrate Eq.~\eqref{eq:dthsm}--\eqref{eq:ptheta} with Matlab's {\tt dde23} integration routine. We then continue such periodic solutions in a parameter with the software DDE-BIFTOOL~\cite{NewDDEBiftool}. This allows us to compute branches of periodic solutions, regardless of whether they are stable or not, to identify their bifurcations, and to continue the respective bifurcation curves in a two-dimensional parameter space. 

\subsection{Excitable smooth theta neuron}

\begin{figure}[t!]
\begin{center}
\includegraphics[width=12cm]{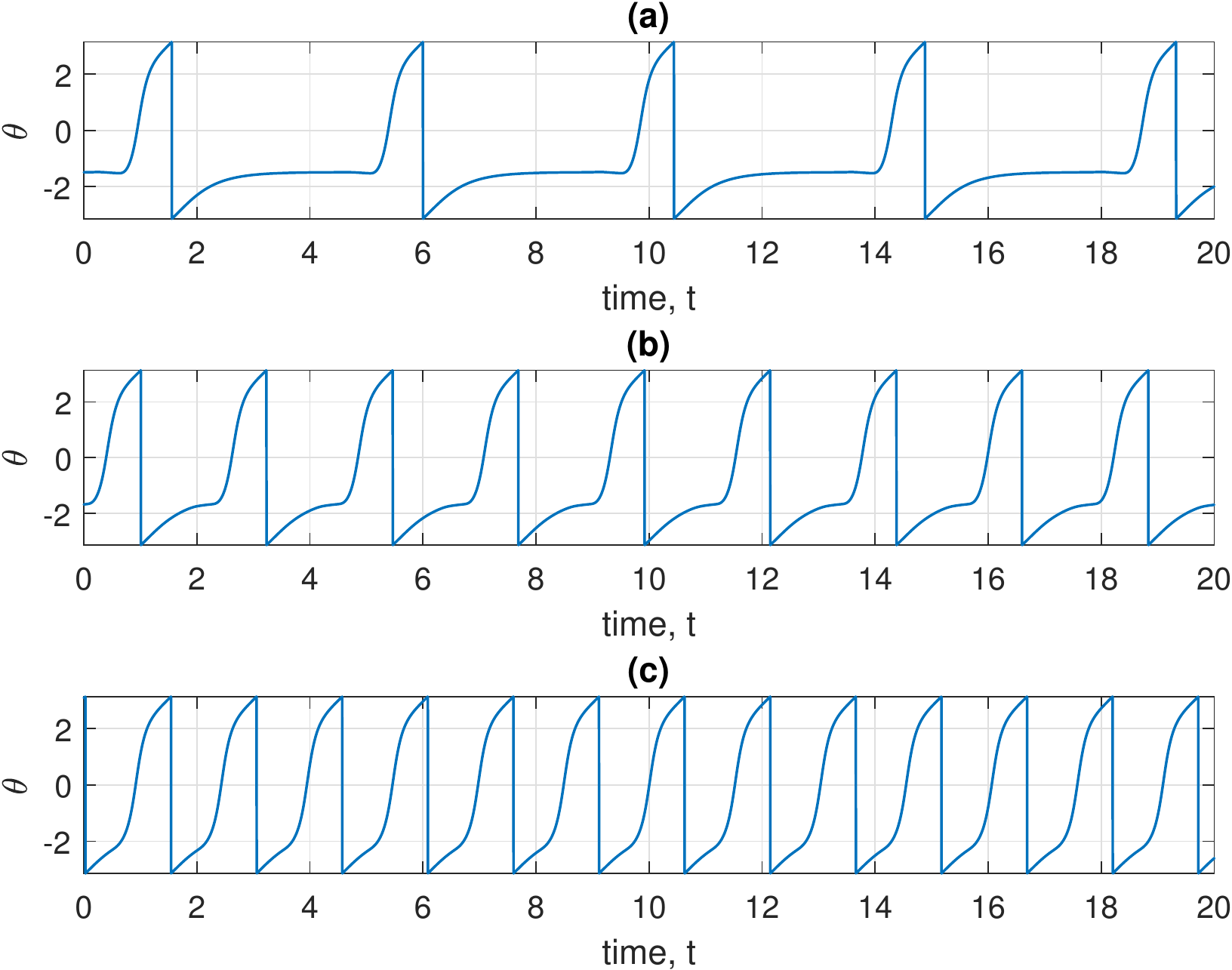}
\caption{Coexisting stable periodic solutions of Eqs.~\eqref{eq:dthsm}--\eqref{eq:ptheta} for $I=-1$, 
$\kappa=2$ and $\tau=4$. Note that when $\theta$ reaches $\pi$ from below it is reset to $-\pi$. Compare with Fig.~\ref{fig:exam}.}
\label{fig:examsm}
\end{center}
\end{figure}

We again first consider the case of an excitable smooth theta neuron with $I<0$, which means that we must have $\kappa>0$ in~\eqref{eq:ptheta} to obtain self-sustained solutions. Figure~\ref{fig:examsm} shows typical coexisting stable periodic solutions of the smooth theta neuron~\eqref{eq:dthsm}--\eqref{eq:ptheta}. Comparison with those in Fig.~\ref{fig:exam} reveals that they have similar shapes, with 1, 2 and 3 spikes per delay interval, respectively. However, the influence of the self-feedback is now smooth, rather than acting instantaneously. Note, in particular, that each periodic solution is smooth across the discontinuty in $\theta$, that is, when $\theta$ is seen as a point on the unit circle. 

\begin{figure}
\begin{center}
\includegraphics[width=9cm]{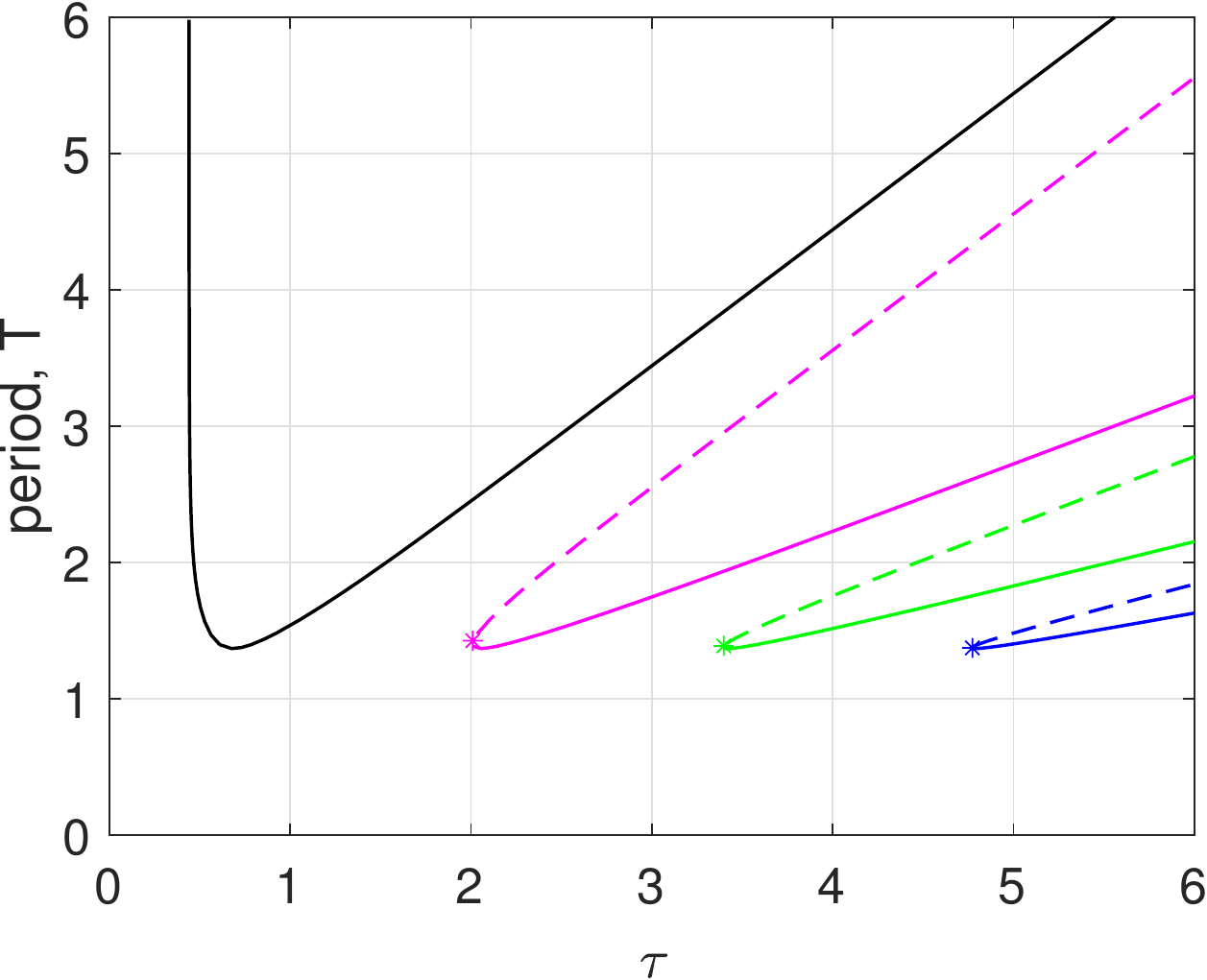}
\caption{Branches of periodic solutions for $n=0,1,2,3$ (left to right) of Eqs.~\eqref{eq:dthsm}--\eqref{eq:ptheta} with negative current $I$, found by numerical continuation and shown by their period $T$ as a function of the delay $\tau$. Periodic solutions are stable along solid curves and unstable along dashed curves; stars indicate saddle-node bifurcations. Here, $I = -1$ and $\kappa=2$, and the three stable solutions in Fig.~\ref{fig:examsm} are those at $\tau=4$. Compare with Fig.~\ref{fig:pertauA}.}
\label{fig:pertausmB}
\end{center}
\end{figure}

\begin{figure}
\begin{center}
\includegraphics[width=10cm]{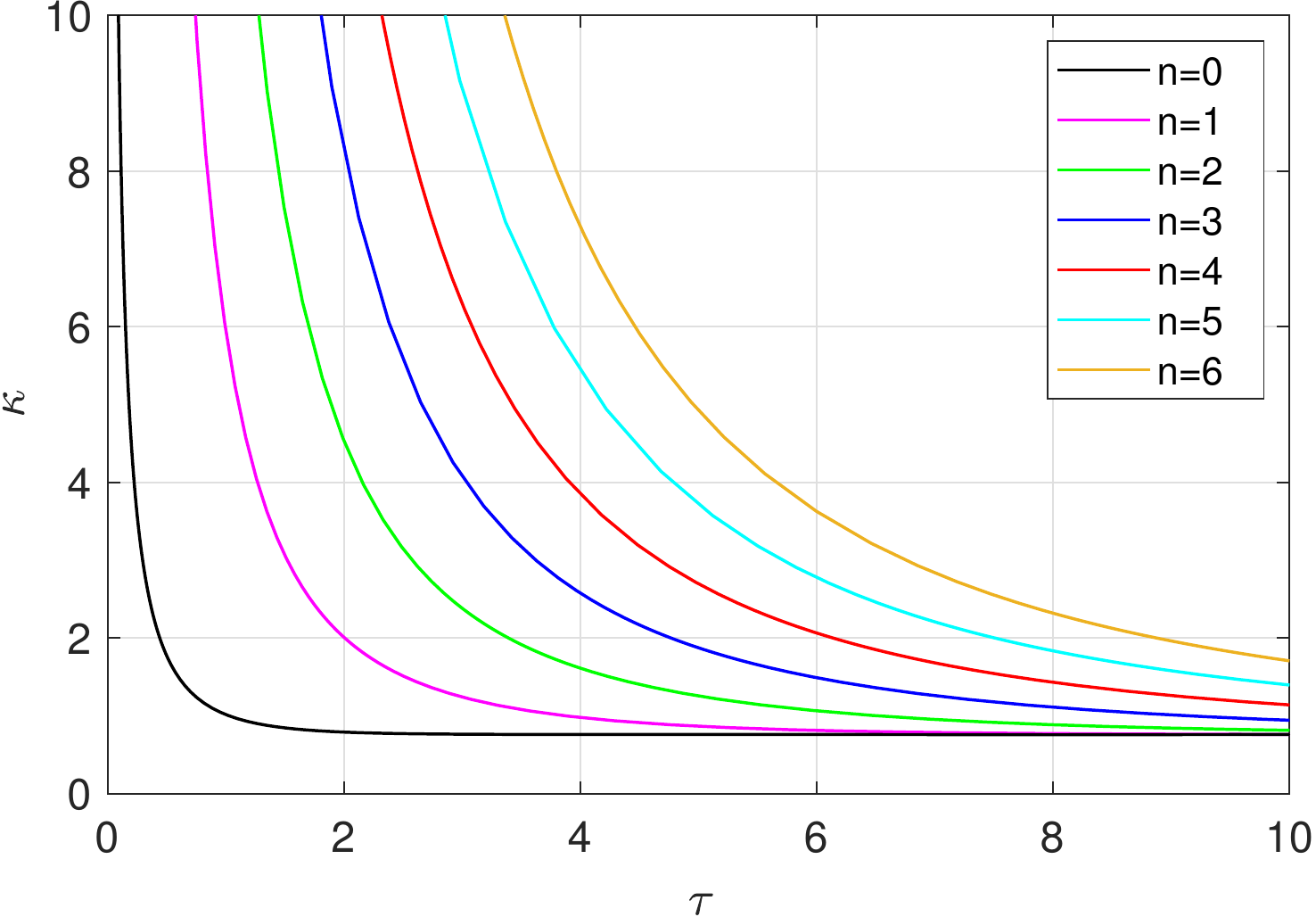}
\caption{Bifurcation curves in the $(\tau,\kappa)$-plane for $n=0,1,2,3,4,5,6$ (see legend) of Eqs.~\eqref{eq:dthsm}--\eqref{eq:ptheta} with negative current $I = -1$. The curve for $n=0$ (black) is the locus of homoclinic bifurcation of the primary branch; it was found via the continuation of a periodic orbit of high period. The curves for $n=1,2,3,4,5,6$ are the saddle-node bifurcations of the secondary branches; they were found by numerical continuation from the respective saddle-node bifurcations in Fig.~\ref{fig:pertausmB}. Compare with Fig.~\ref{fig:kaptauA}.}
\label{fig:snsmex}
\end{center}
\end{figure}

Figure~\ref{fig:pertausmB} shows the first four branches of the periodic solutions obtained by numerical continuation as $\tau$ is varied; observe that for $\tau = 4$ there are indeed three simultaneously stable solutions as shown in Fig.~\ref{fig:examsm}. The branches of periodic orbits are distinct and, as before, can be indexed in Fig.~\ref{fig:pertausmB} based on the number of times $\theta$ increases through $\pi$ in one delay period; this corresponds to the number $n$ of additional spikes per delay interval. We observe qualitatively the same behaviour as in Fig.~\ref{fig:pertauA} for instantaneous self-feedback: the single-spike basic branch for $n=0$ is a function over $\theta$, while for $n \geq 1$ there is a stable and an unstable branch that meet at a saddle-node bifurcation. Figure~\ref{fig:snsmex} shows the bifurcations for $n=0,1,2,3,4,5,6$ as curves in the $(\tau,\kappa)$-plane, which must be found by numerical continuation. The primary branch approaches a homoclinic bifurcation, which was identified, and then continued, as a single-spike orbit of sufficiently large period. This resulted in the curve for $n=0$ of Fig.~\ref{fig:snsmex}, to the left and below of which no periodic solutions exist. The curves for $n=1,2,3,4,5,6$ are curves of saddle-node bifurcations, and they have been continued from the identified points of saddle-node bifurcations on the secondary branches in Fig.~\ref{fig:pertausmB}. We observe that the bifurcation curves in Fig.~\ref{fig:snsmex} are also in perfect qualitative agreement with the corresponding analytically obtained curves in Fig.~\ref{fig:kaptauA} for instantaneous self-feedback.

We conclude that system~\eqref{eq:dth} with Dirac delta function predicts correctly the observed behaviour of the smooth DDE model~\eqref{eq:dthsm}--\eqref{eq:ptheta} for the case that the theta neuron is excitable. This is quite remarkable given that the feedback spike $P(\theta)$ in~\eqref{eq:ptheta} is of a considerable width and not close to a Dirac delta function. 

\begin{figure}
\begin{center}
\includegraphics[width=11.5cm]{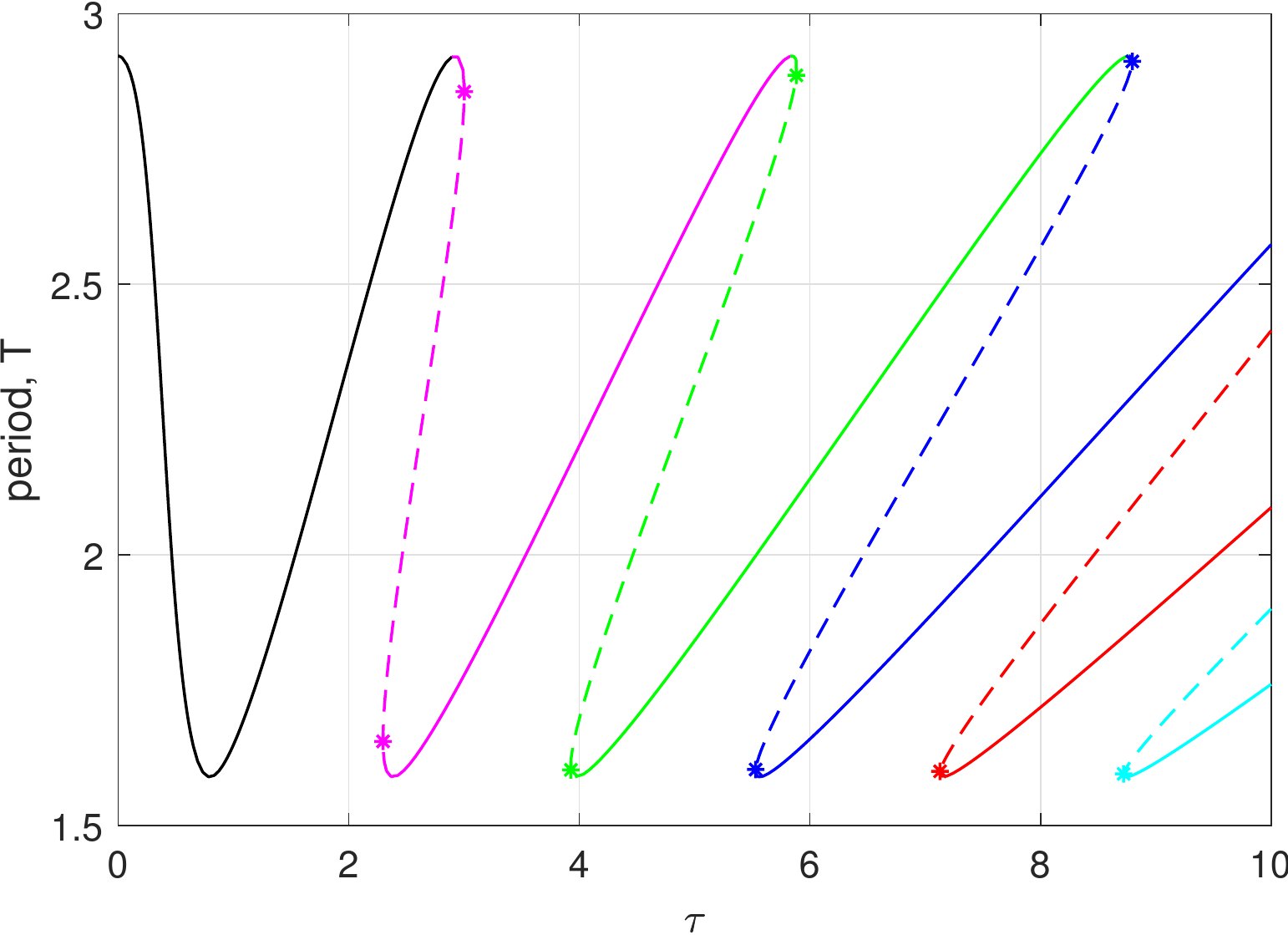}
\caption{Periodic solutions of~\eqref{eq:dthsm}-\eqref{eq:ptheta} with positive current $I$ and positive $\kappa$, shown by period $T$ as a function of the delay $\tau$. This single branch was found by numerical continuation and has been split up at its maxima into segments associated with additional firing events $n=0,1,2,3,4,5$. Periodic solutions are stable along solid curves and unstable along dashed curves; stars indicate saddle-node bifurcations. Here, $I = 1$ and $\kappa=1$. Compare with Fig.~\ref{fig:pertauposA}.}
\label{fig:pertauposB}
\end{center}
\end{figure}

\begin{figure}
\begin{center}
\includegraphics[width=11cm]{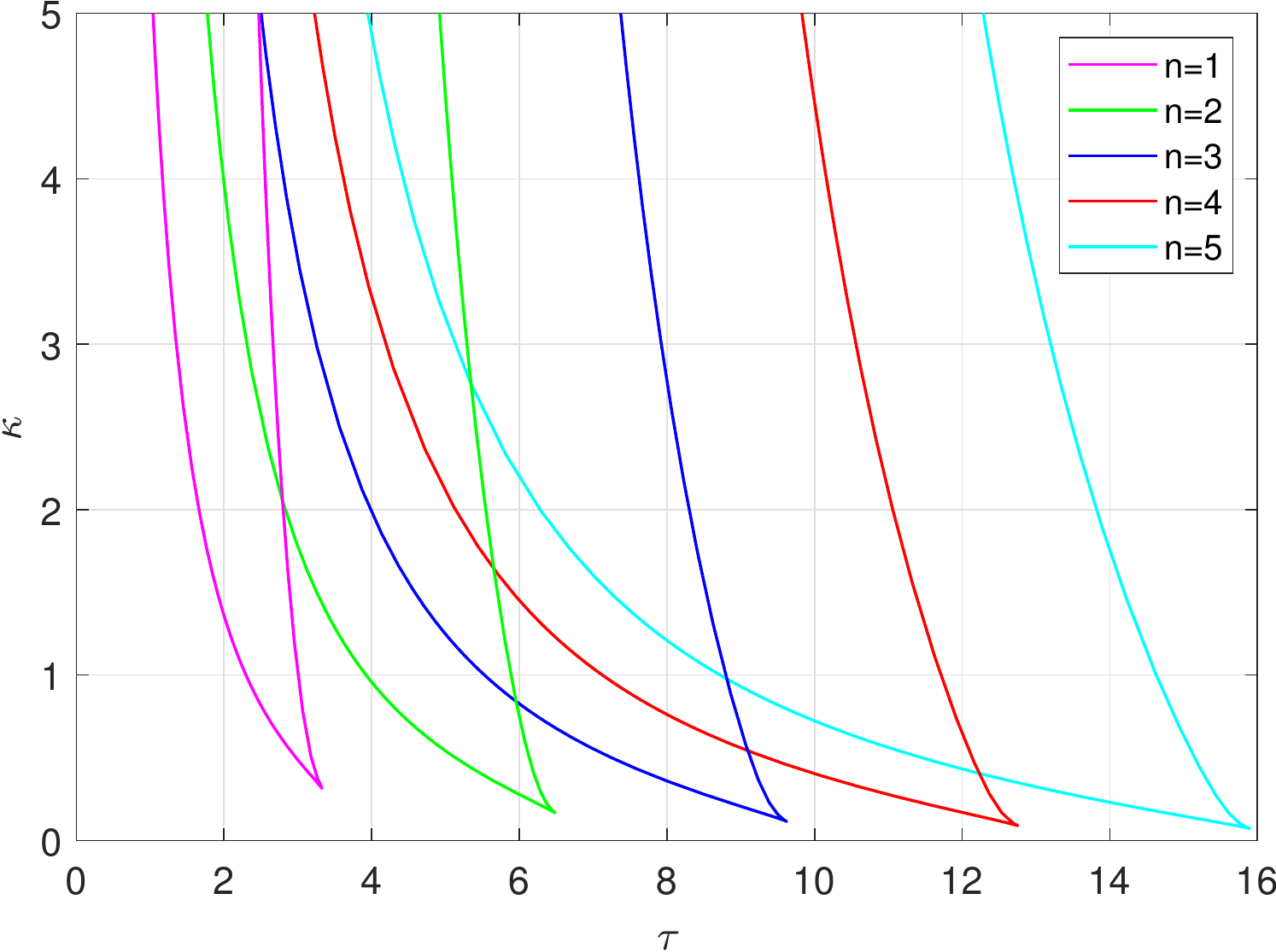}
\caption{Saddle-node bifurcation curves in the $(\tau,\kappa)$-plane with $n=1,2,3,4,5$ (see legend) of Eqs.~\eqref{eq:dthsm}--\eqref{eq:ptheta} with positive current $I = 1$, found by numerical continuation from the respective saddle-node bifurcations in Fig.~\ref{fig:pertauposB}. Compare with Fig.~\ref{fig:kaptaupos} for $\kappa > 0$.}
\label{fig:snsmin}
\end{center}
\end{figure}

\subsection{Intrinsically oscillating smooth theta neuron}
\label{sec:smoopos}

We now consider the smooth theta neuron~\eqref{eq:dthsm}--\eqref{eq:ptheta} with positive $I$ when, in the absence of coupling, the neuron fires periodically with period $\pi/\sqrt{I}$. When $\kappa$ is positive the self-coupling is excitatory and numerical continuation started from a stable periodic solution results in the single branch of periodic solutions shown in Fig.~\ref{fig:pertauposB}. While the periodic solution varies smoothly along the branch, successive segments of it can still be associated with an increasing number of firings within a single delay period. We find that the transitions between consecutive numbers of firings take place at the maxima, and we distinguish the respective segments again by colour in Fig.~\ref{fig:pertauposB} to represent the associated value of additional spikes $n=0,1,2,3,4,5$. This highlights the very good agreement with Fig.~\ref{fig:pertauposA} for~\eqref{eq:dth} with Dirac delta function. Excellent qualitative agreement is also observed when we continue the saddle-node bifurcations identified in Fig.~\ref{fig:pertauposB}. The resulting bifurcation curves are shown in Fig.~\ref{fig:snsmin} for $n=1,2,3,4,5$, and should be compared with the top half of Fig.~\ref{fig:kaptaupos}.

\begin{figure}
\begin{center}
\includegraphics[width=11.5cm]{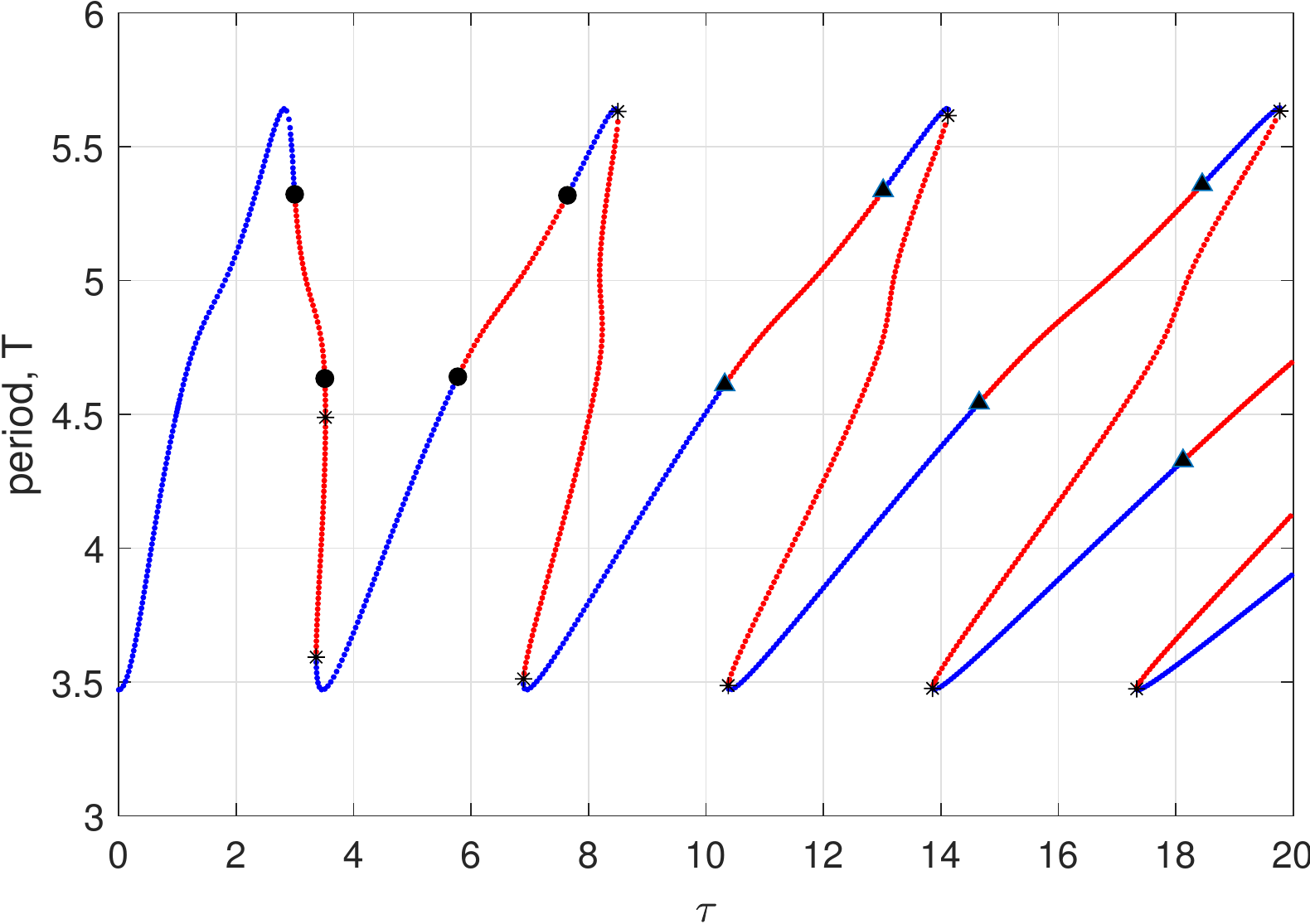}
\caption{Periodic solutions of~\eqref{eq:dthsm}-\eqref{eq:ptheta} with positive current $I$ and negative $\kappa$, shown by period $T$ as a function of the delay $\tau$. This single branch was found by numerical continuation while detecting saddle-node bifurcations (stars), 
period-doubling bifurcations (dots), and Neimark-Sacker bifurcations (triangles); solutions along the branch are stable along blue segments and unstable along red segements. Here, $I = 1$ and $\kappa=-1$. Compare with Fig.~\ref{fig:pertaunegA}.}
\label{fig:pertaunegB}
\end{center}
\end{figure}

\begin{figure}
\begin{center}
\includegraphics[width=12cm]{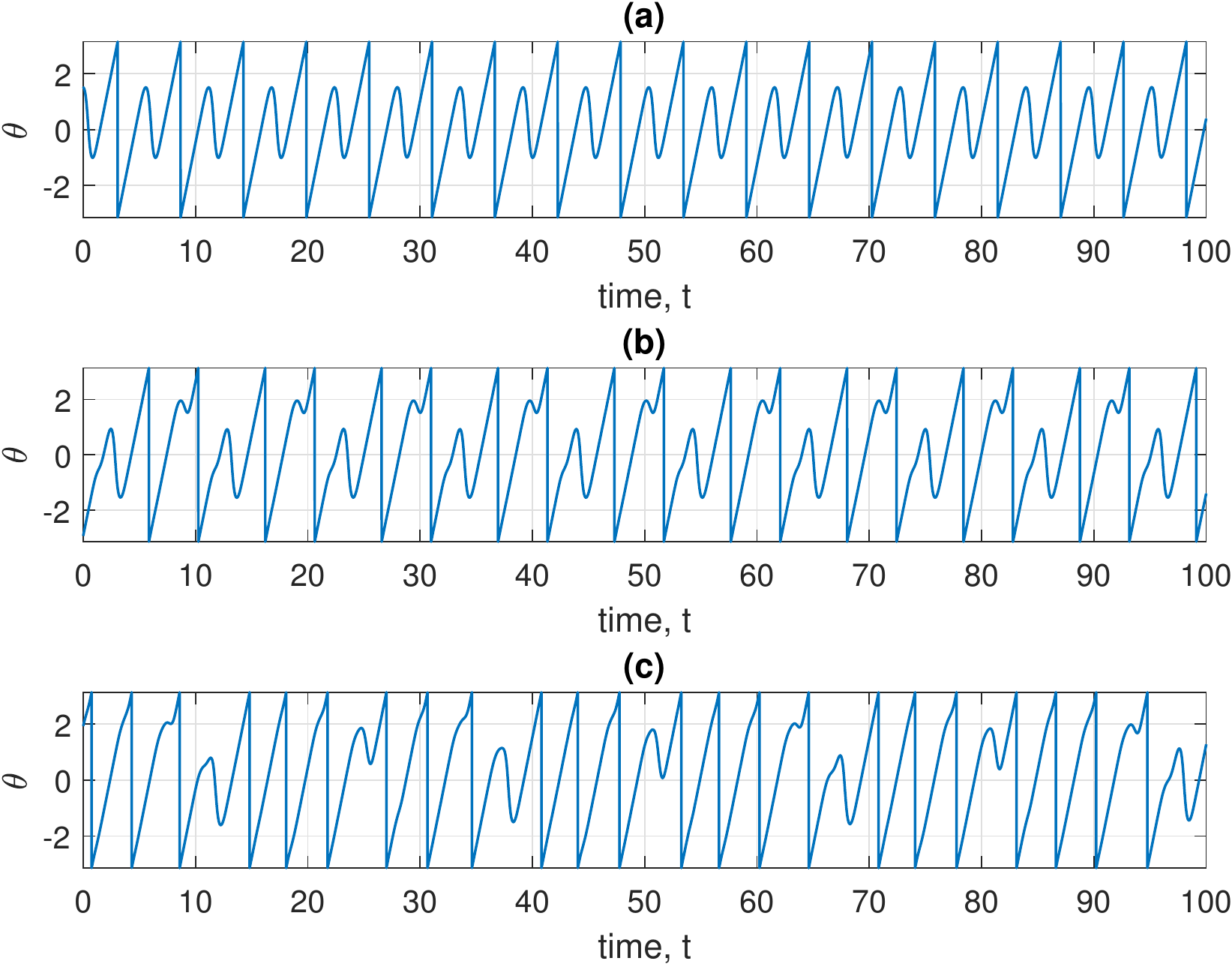}
\caption{Solutions of~\eqref{eq:dthsm}--\eqref{eq:ptheta} with $I = 1$ and $\kappa=-1$ for $\tau=2.9$ (a), $\tau=3.05$ (b) and $\tau=3.3$ (c) illustrate period-doubling to chaotic spiking.}
\label{fig:chaos}
\end{center}
\end{figure}

We finally consider the smooth theta neuron~\eqref{eq:dthsm}--\eqref{eq:ptheta} with positive $I$ for the case $\kappa < 0$ when the self-coupling is inhibitory. As was expected from the equivalent case of the theta neuron~\eqref{eq:dth} with Dirac delta function in Fig.~\ref{fig:pertaunegA}, there is again a single branch of periodic solutions, which is shown in Fig.~\ref{fig:pertaunegB}. More precisely, on this branch one finds orbits with a different number $n$ of additional spikes (this is not shown by colour here), and these all feature a pair of saddle-node bifurcations. We remark that the corresponding curves of saddle-node bifurcations in the $(\kappa,\tau)$-plane (not reproduced here) are qualitatively as those in the bottom half of Fig.~\ref{fig:kaptaupos}.
However, additional points of bifurcation are now also detected during the continuation of the branch in Fig.~\ref{fig:pertaunegB}, specifically, period-doubling and Neimark-Sacker bifurcations. This leads to additional intervals, between pairs of these two types of bifurcations, where the respective periodic orbit is unstable. In particular, there are now ranges of $\tau$ for which neither the primary orbit nor the secondary periodic orbit is stable. This opens up the possibility of more complicated periodic and even chaotic dynamics. This is demonstrated in Fig.~\ref{fig:chaos} with three time series near the first period-doubling bifurcation at $\tau \approx 2.99$, showing periodic, period-doubled and chaotic spiking, respectively. The statement that the solution in panel (c) is chaotic has been checked by verifying that it features a positive Lyapunov exponent. 

Overall, we find a more diverse picture for the case that the theta neuron is intrinsically oscillating. For excitatory self-feedback there is still excellent agreement between the theta neuron~\eqref{eq:dth} with Dirac delta function and the smooth theta neuron~\eqref{eq:dthsm}--\eqref{eq:ptheta}. For inhibitory self-coupling, on the other hand, we still find a single branch with pairs of saddle-node bifurcations as for system~\eqref{eq:dth}, but the smooth DDE~\eqref{eq:dthsm}--\eqref{eq:ptheta} now features additional bifurcations that lead to more complicated dynamics in certain ranges of the delay $\tau$, including chaotic spiking. It might be argued that this more complex behaviour results from the somewhat non-physical nature of the model in this parameter regime of positive current $I$ and negative self-feedback $\kappa$, which is why we do not investigate it further here. 

\section{Discussion}
\label{sec:disc}

We have studied the dynamics of a single theta neuron with delayed self-coupling in the form of a delta function. Because the dynamics can be solved explicitly between times at which the feedback occurs, we were able to analytically describe periodic orbits and determine their stability in terms of the roots of a finite-order polynomial. The only possible bifurcations of periodic orbits are saddle-node bifurcations, and we gave explicit expressions for curves in two-dimensional parameter space along which they are found. In this way, we provided a complete description of the types of spiking solutions, where they occur and their stability. 

The theta neuron with delayed delta-function self-coupling can be thought of as a ``normal form'' of an excitable system subject to delayed self-feedback. By this we mean that the system is solvable explicitly, while still capturing the essentials of the behaviour of other excitable systems with non-instantaneous feedback that do not have analytical solutions. Indeed, the kind of dynamics and bifurcation structures presented here have been found (by means of numerical techniques) in different contexts as well~\cite{wedslo21,GarbinNC15}. To test the predictive power of our results for delta-function feedback, we investigated a single theta neuron with smooth delayed self-coupling. This showed that for excitatory feedback the dynamics are qualitatively the same for both types of intrinsic dynamics (excitable and intrinsically oscillating). In the intrinsically oscillating regime with inhibitory feedback the basic structure of periodic solutions is still as predicted, but we found additional, more complex behaviour generated by period-doubling and Neimark-Sacker bifurcations. These bifurcations have also been found in \cite{TerrienPRE21} on branches of periodic spiking in the context of pulse-timing symmetry breaking in a nanolaser with optical feedback.

We now discuss related work that concerns networks of neurons. Several groups have studied infinite networks of QIF (or equivalently, theta) neurons with delayed feedback~\cite{devmon18,pazmon16,ratpyr18}. Devalle et al.~\cite{devmon18} considered an infinite network of excitable QIF neurons with delayed delta function coupling. The synchronous state in that network is described by the dynamics of one neuron, as we study here. These authors considered the case of a single spike in the delay interval and derived an expression equivalent to~\eqref{eq:hom} giving the minimum value of the feedback strength for which a synchronous periodic solution can exist. They also analysed the stability of such a solution but obtained different results from us, as there are instabilities in infinite networks that cannot occur for a single self-coupled neuron.  In similar work, Paz{\'o} and Montbri{\'o}~\cite{pazmon16} considered the same network but in the intrinsically firing regime. These authors derived an expression describing the existence of a synchronous state equivalent to~\eqref{eq:existA}. They analysed the stability of such a solution and again obtained different results from us, due to the instabilities mentioned above. In relation to the occurence of more complex spiking behaviour dynamics, chaotic dynamics was found, also in ~\cite{pazmon16}, in an intrinsically oscillating infinite network of identical QIF neurons with delayed inhibitory delta function feedback. However, this chaotic behaviour required that the neurons are not synchronised, meaning that it is not equivalent to the dynamics of the single-neuron model we studied here. Chaotic dynamics were also found in a network of three intrinsically oscillating theta neurons with nondelayed smooth inhibitory feedback~\cite{lai18A}, but this is due to the reversibility of the network's dynamics.

Possible generalisations of our work include the study of two coupled neurons~\cite{dahhil09} or a ring of unidirectionally coupled neurons~\cite{kliluc18}. For more complex networks, the fact that we can explicitly solve for the dynamics of an uncoupled neuron means that we could efficiently implement event-based simulations: jumping straight from one firing event to the next without having to numerically integrate differential equations between them~\cite{bre06}. Such schemes are very efficient and have been implemented in~\cite{monwol10,ernpaw95}, for example. Another possibility is to consider an infinite network of heterogeneous theta neurons for which the Ott/Antonsen ansatz~\cite{lai14A,ottant08} could be used to derive a single complex delay differential equation for the network's order parameter, as has been done for delayed Kuramoto oscillators~\cite{leeott09,lai16A}.







\bigskip

{\bf Acknowledgements:}
The authors thank Stefan Ruschel for helpful conversations about delay differential equations.




\vspace*{5mm}
{\bf Affiliations:}

\vspace*{2mm}

\end{document}